\documentclass[journal]{IEEEtran}

\usepackage{amsmath}
\usepackage{amsfonts}
\usepackage{amssymb}
\usepackage{graphicx}
\usepackage{enumitem}
\usepackage{algorithm}
\usepackage{algpseudocode}
\usepackage{array}
\usepackage{cite}
\usepackage{color}
\newtheorem{defn}{Definition}
\newtheorem{lem}{Lemma}

\newtheorem{ass}{Assumption}
\newtheorem{thm}{Theorem}

\newtheorem{rmk}{Remark}

\newcommand{\cm}[1]{{\color{black} {#1}}}

\begin{document}
\title{\huge Decentralized Schemes with Overlap for Solving Graph-Structured Optimization Problems}
\author{Sungho~Shin,
  Victor~M.~Zavala,
  and Mihai~Anitescu 
  \thanks{S. Shin is with the Department of Chemical and Biological Engineering, University of Wisconsin-Madison, Madison, WI 53706 USA (e-mail: sungho.shin@wisc.edu).}
  \thanks{V. M. Zavala is with the Department of Chemical and Biological Engineering, University of Wisconsin-Madison, Madison, WI 53706 USA, and also with the Mathematics and Computer Science Division, Argonne National Laboratory, Lemont, IL 60439 USA (e-mail: victor.zavala@wisc.edu).}
  \thanks{M. Anitescu is with the Mathematics and Computer Science Division, Argonne National Laboratory, Lemont, IL 60439 USA, and also with the Department of Statistics, University of Chicago, Chicago, IL 60637 USA (e-mail: anitescu@mcs.anl.gov)}
}
\maketitle

\begin{abstract}
We present a new algorithmic paradigm for the decentralized solution of graph-structured optimization problems that arise in the estimation and control of network systems. A key and novel design concept of the proposed approach is that it uses overlapping subdomains to promote and accelerate convergence. We show that the algorithm converges if the size of the overlap is sufficiently large and that the convergence rate improves exponentially with the size of the overlap. The proposed approach provides a bridge between fully decentralized and centralized architectures and is flexible in that it enables the implementation of asynchronous schemes, handling of constraints, and balancing of computing, communication, and data privacy needs. {The proposed scheme is tested in an estimation problem for a 9241-node power network and we show that it outperforms the alternating direction method of multipliers.}
\end{abstract}

\begin{IEEEkeywords}
  decomposition, parallel computing, asynchronous 
\end{IEEEkeywords}
\IEEEpeerreviewmaketitle

\section{Introduction}
\begin{figure*}[!t]
  \centering
  \includegraphics[width=\textwidth]{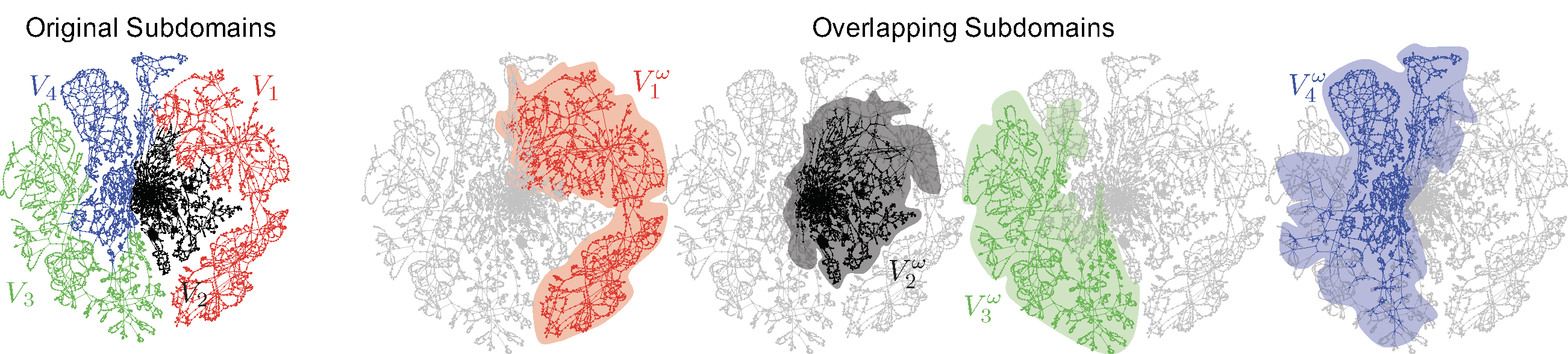}
  \caption{Sketch of the decomposition of a graph into overlapping subdomains.}
  \label{fig:abstract-2}
\end{figure*}

Diverse centralized and decentralized architectures for control and estimation of large-scale networks have been reported in the literature. Centralized architectures achieve high performance due to their ability to capture interconnections among all the network components but their scalability is hindered by computing, communication, and data privacy issues. Decentralized schemes can help mitigate these obstacles. In these schemes, the original network domain is partitioned into multiple tractable subdomains that are coordinated to try to achieve optimal performance. Decomposition schemes such as the alternating direction method of multipliers (ADMM) \cite{boyd2011distributed}, Lagrangian decomposition (price coordination) \cite{giselsson2013accelerated,kozma2013distributed}, Benders decomposition \cite{geoffrion1972generalized}, and Gauss-Seidel schemes \cite{shin2018multi} have been explored in the literature. {ADMM, in particular, has been widely used in network applications \cite{kozma2015benchmarking}.  This approach has proven to be scalable and flexible but it is well-known that its convergence can be rather slow and that it is difficult to tune \cite[Section 3.2.2]{boyd2011distributed}.}

A powerful distributed computing paradigm that has been widely studied in the context of partial differential equations (PDEs) is {\em overlapping domain decomposition} (often referred to as overlapping Schwarz)  \cite{mathew2008domain,quarteroni1999domain,toselli2006domain,cai1992domain,dryja1994domain}. Convergence can be accelerated with this scheme by increasing the size of the overlap  \cite{dryja1987additive,mathew2008domain,schwarz1870ueber}.  In this work, we establish a Schwarz-like, decentralized computing architecture for solving optimization problems that arise in the control and estimation of network systems {(we call these {\em graph-structured optimization problems})}. With the proposed approach, the original underlying graph domain is partitioned into multiple subdomains with a given degree of overlap (see Fig. \ref{fig:abstract-2}). At each coordination step, the subproblems associated with the subdomains are solved in parallel by using current information of the neighbors in the overlapping subdomains. Our analysis shows that convergence is guaranteed if the size of the overlap is sufficiently large and that the convergence rate improves exponentially with the size of the overlapping region. Moreover, we show that convergence can be guaranteed under asynchronous coordination. The proposed architecture provides a bridge between fully decentralized schemes (with no overlapping region) and fully centralized schemes (the overlapping region is the entire domain). Moreover, the architecture provides design flexibility to balance computing, communication, and privacy needs. {We provide numerical results that illustrate theoretical properties and compare the performance of the scheme with that of ADMM. The results suggest that the proposed scheme is efficient at solving large-scale problems and outperforms ADMM.}

\cm{Classical work on overlapping Schwarz has focused on linear algebra systems that arise from the discretization of PDEs on regular meshes and more recent work has focused on generalizations. In \cite{chan1994additive}, for instance, convergence is analyzed for unstructured meshes by treating the operator matrix of the unstructured mesh as an approximate operator for the regular mesh. Schwarz schemes for general linear algebra systems have also been proposed and these are available in powerful scientific computing packages such as {\tt PETSc} \cite{balay2019petsc}. Most work reported on general linear systems, however, has only analyzed empirical convergence behavior   \cite{cai1996overlapping,cai1999restricted,chan1994domain,gander2014block}. The work in  \cite{frommer1999weighted,jones1996two} provides theoretical convergence analysis for general linear systems but does not establish a relationship between the convergence rate and the size of the overlap. In this context, we highlight that the technical contributions of this work are: (i) We provide a theoretical convergence analysis for a Schwarz overlapping scheme for general (graph-structured) positive definite linear systems under synchronous and asynchronous settings. (ii) We establish an explicit dependence between the convergence rate and the size of the overlap. (iii) We analyze the Schwarz scheme from an optimization perspective. This highlights applications of Schwarz schemes in control and estimation of network systems and points towards extensions that can handle constraints and nonlinearities. To the best of our knowledge,  the use of Schwarz schemes for solving optimization problems has not been reported in the literature.} 

The paper is structured as follows. In Section \ref{sec:setting}, we introduce the proposed setting and basic notation. In Section \ref{sec:dom}, we introduce the notion of the graph-structured matrix and derive basic analytical tools that help analyze the proposed architecture. Convergence analysis under synchronous and asynchronous coordination schemes is presented in Section \ref{sec:dec}. In Section \ref{sec:impl}, we discuss implementation details of the architecture in a distributed computing environment. In Section \ref{sec:cstudy}, we demonstrate the performance of the proposed algorithm in a large-scale state estimation problem. Conclusions and directions of future work are presented in Section \ref{sec:con}.

\section{Basic Notation and Setting}\label{sec:setting}

The set of real numbers and the set of complex numbers are denoted by $\mathbb{R}$ and $\mathbb{C}$, respectively. We define $\mathbb{R}_{>0}:=(0,+\infty)$, $\mathbb{R}_{\geq 0}:=[0,+\infty)$. The absolute values of real numbers, the magnitudes of complex numbers, and the cardinality of sets are denoted by $|\cdot|$. The ceiling operator is denoted by $\lceil\cdot\rceil$. The $i$th component of a vector and the $(i,j)$th component of a matrix are denoted by $(\cdot)_i$ and  $(\cdot)_{i,j}$, respectively. The $j$th column vector and the $i$th row vector of a matrix are denoted by $(\cdot)_{:,j}$ and $(\cdot)_{i,:}$, respectively. The transpose of a matrix or a vector is denoted by $(\cdot)^T$. We use the notation $(x_1,x_2,\cdots,x_n) = \begin{bmatrix}x_1^T & x_2^T &\cdots &x_n^T\end{bmatrix}^T$, where $x_1,x_2\cdots,x_n$ are scalars or column vectors, and $\{x_i\}_{i\in I} := (x_{i_1},x_{i_2},\cdots,x_{i_n})$, where $I=\{i_1<i_2<\cdots<i_n\}$. We use the following syntax:
\begin{align}
  \{x_{ij}\}_{i\in I,j\in J}:=
  \begin{bmatrix}
    x_{i_1 j_1}&\cdots &x_{i_1 j_m}\\
    \vdots&\ddots &\vdots\\
    x_{i_n j_1}&\cdots &x_{i_n j_m},            
  \end{bmatrix}
\end{align}
where $I=\{i_1<\cdots<i_n\}$ and $J=\{j_1<\cdots<j_m\}$. Furthermore, $\{x_{i}\}_{i\in I,:}:=\begin{bmatrix}x_{i_1}^T& \cdots & x_{i_n}^T\end{bmatrix}^T$, where $x_i$ are row vectors and $\{x_{j}\}_{:,j\in J}:=\begin{bmatrix}x_{j_1}& \cdots & x_{j_n}\end{bmatrix}$, where $x_j$ are column vectors. For simplicity, $\{(\cdot)_i\}_{i\in I}$ is replaced by $\{\cdot\}_{I}$. The spectral radius  of a matrix (the magnitude of the eigenvalue with the largest magnitude) is denoted by $\rho(\cdot)$. Infinity norms of vectors and induced infinity norms of matrices are denoted by $\Vert \cdot \Vert_\infty$, where the induced infinite matrix norm of $A\in\mathbb{R}^{n\times n}$ is defined by $\Vert A\Vert_\infty :=\max_{\Vert x \Vert_\infty=1} \Vert Ax \Vert_\infty$, or equivalently (see \cite[Proposition 5.6.5]{horn1990matrix}), $\Vert A \Vert_\infty := \max_{1\leq i\leq n} \sum_{j=1}^n |(A)_{i,j}|$. The set of neighbors of node $i\in{V}$ on an undirected graph $G({V},{E})$ is denoted by $\mathcal{N}_G(i):=\{j\in{V}\mid \{i,j\}\in{E}\}$. The distance between vertices $i,j\in{V}$ on a graph $G(V,E)$---the smallest number of edges on a path between the two nodes ---is denoted by $d_G(i,j)$. We extend this concept to the distance between sets $X,Y\subseteq{V}$ of vertices, which is defined as $d_G(X,Y) := \min \{d_G(x,y)\mid x\in X, y\in Y\}$ (when $X$ or $Y$ are single nodes, we consider them as singletons).

We introduce a motivating graph-structured optimization problem that captures various applications such as optimal power flow, power system state estimation, and PDE control. Specifically, we consider a problem on a graph $G({V},{E})$ of the form:
\begin{subequations}\label{eqn:orig-graph}
  \begin{align}
    \min_{x,u,v}\ &\sum_{i\in V} \left(\frac{1}{2}q_i x_i^2 -f_i x_i + \frac{1}{2}r_i u_i^2 \right) + \sum_{\{i,j\}\in E} \frac{1}{2}s_{ij} v_{ij}^2\\
    \text{s.t.}\ &a_{ii}x_i-\sum_{j\in\mathcal{N}_G(i)} a_{ij}x_j = u_i \quad i\in{V}\\
    & b_{ij}(x_i-x_j) = v_{ij} \quad \{i,j\}\in{E}.    \label{eqn:direction}
  \end{align}
\end{subequations}
Here $x_i\in\mathbb{R}$ is the state of node $i$; $u_i\in\mathbb{R}$ is the input of node $i$; $v_{ij}\in\mathbb{R}$ is the state of edge $\{i,j\}\in E$; and $q_i,r_i\in\mathbb{R}_{\geq 0}$ and $f_i\in\mathbb{R}$ for $i\in V$ and $s_{ij}\in\mathbb{R}_{\geq 0}$ and $a_{ij},b_{ij}\in\mathbb{R}$ for $\{i,j\}\in E$ are either data, model parameters, or the components in the objective functions. For convenience, we assume that a direction is assigned to each edge so that $s_{ij}$, $b_{ij}$, and $v_{ij}$ for edge $\{i,j\}\in E$ can be uniquely defined. We label the vertices by $V:=\{1,2,\cdots,n\}$.  Such a problem can be equivalently written in vector form as:
\begin{subequations}\label{eqn:orig}
  \begin{align}
    \min_{x,u}\,&\frac{1}{2} x^T Q x - f^T x +\frac{1}{2}u^T Ru  + \frac{1}{2}v^TS v \label{eqn:orig-obj}\\
    \text{s.t.}\,&Ax = u,\quad Bx = v \label{eqn:orig-con},
  \end{align}
\end{subequations}
where $Q,R,A\in\mathbb{R}^{n\times n}$, $B\in\mathbb{R}^{m\times n}$, $S \in \mathbb{R}^{m \times m}$, $x,u\in\mathbb{R}^n$, and $v\in\mathbb{R}^m$.  Upon elimination of the equality constraints \eqref{eqn:orig-con}, one can derive the equivalent quadratic program (QP):
\begin{align}\label{eqn:elim}
  \min_{x}\,&\frac{1}{2}x^T H x -f^T x.
\end{align}
where we denote $H:=Q+A^T R A+ B^T S B$, and we assume that $H$ is positive definite (PD). Consequently, the solution of \eqref{eqn:elim} exists, is unique, and can be obtained by solving:
\begin{align}\label{eqn:lin}
  H x = f.
\end{align}
An important observation is that one can derive solution schemes that equivalently operate in the linear algebra space  \eqref{eqn:lin}, in the reduced QP space \eqref{eqn:elim}, or in the full QP space \eqref{eqn:orig}.  Operating in the full QP space facilitates implementation because this requires less-intrusive manipulations but convergence analysis is in general easier in the linear algebra space. 

The graph structure is embedded in matrices $A$ and $B$. In particular, $(A)_{i,j}$ is nonzero only if $d_G(i,j)\leq 1$ and there exists $k$ such that $(B)_{k,i},(B)_{k,j}\neq 0$ only if $\{i,j\}\in E$. This implies that the sparsity of $A$ and $B$ can be characterized by the structure of the graph $G$. The sparsity structure is partially preserved in $H$; specifically, if $(H)_{i,j}$ is nonzero, then $d_G(i,j)\leq 2$. The sparsity induced by the graph is the key property that will be exploited in the analysis of the proposed decentralized architecture.

\section{Graph-Induced Matrix Properties}\label{sec:dom}

In this section, we introduce the concept of the {\em generalized matrix bandwidth}, which is an essential tool for analyzing the convergence properties of the proposed decentralized architecture. The conventional matrix bandwidth is defined as the maximum range of a diagonally bordered band such that the components beyond the band are all zeros. Here, we define the generalized matrix bandwidth by finding the maximum range of a band {\em induced by a graph} such that the components beyond that band are all zeros. We will also see that this notion can be interpreted as a diffusion process that captures long-range coupling between different elements of the network. 

\begin{defn}\label{defn:bandwidth}
  Suppose that we have a matrix $A\in\mathbb{C}^{n\times n}$ and a graph $G(V,E)$ with $\{1,2,\cdots,n\}\subseteq V$.
  \begin{enumerate}[label=(\alph*)]
  \item A bandwidth of $A$ with respect to $G$ is the smallest nonnegative integer $\mathcal{B}_G(A)$ such that $(A)_{i,j}= 0$ for any $i,j\in\{1,2,\cdots,n\}$ and $d_G(i,j)>\mathcal{B}_G(A)$.
  \item {$A$ is graph-structured with respect to $G$ if the bandwidth $\mathcal{B}_G(A)$ of $A$ with respect to $G$ is sufficiently smaller than the diameter of $G$.}
  \end{enumerate}
\end{defn}
{Here, the diameter of $G$ is defined as the maximum distance between two nodes in $V$.}
The generalized bandwidth changes with the choice of graph $G$; for example, when a linear graph is chosen for $G$, the bandwidth with respect to $G$ reduces to the standard definition of the matrix bandwidth. On the other hand, the bandwidth is always less than or equal to one if a complete graph (a graph in which each pair of graph vertices is connected by an edge) is chosen for $G$. 

{Graph-induced structures are preserved under addition and multiplication. In particular, the generalized bandwidth does not abruptly increase under matrix summations and multiplications. This observation is rigorously stated as follows.}
\begin{lem}\label{lem:bandwidth}
  The following holds for $A,B\in\mathbb{C}^{n\times n}$ and a graph $G(V,E)$ with $\{1,2,\cdots,n\}\subseteq V$.
  \begin{subequations}
    \begin{align}
      \mathcal{B}_G(A+B)&\leq \max\left(\mathcal{B}_G(A),\mathcal{B}_G(B)\right)\label{eqn:bandwidth-a}\\
      \mathcal{B}_G(AB)&\leq \mathcal{B}_G(A)+\mathcal{B}_G(B)\label{eqn:bandwidth-b}
    \end{align}
  \end{subequations}
\end{lem}
The proof is given in Appendix \ref{prop:bandiwdth-pf}.

The generalized bandwidth (Definition \ref{defn:bandwidth}) and its property (Lemma \ref{lem:bandwidth}) can be used to find componentwise bounds for the inverse of graph-structured matrices. If $H$ is PD, we can establish the following result. 

\begin{thm}\label{thm:domdec}
  Consider a PD matrix $H\in\mathbb{R}^{n\times n}$ and a graph $G(V,E)$ with $\{1,2,\cdots,n\}\subseteq V$. Suppose that any eigenvalue $\lambda$ of $H$ satisfies $\lambda\in[\lambda_{\min},\lambda_{\max}]$ for some $\lambda_{\min},\lambda_{\max}\in\mathbb{R}_{>0}$. Then the following holds for any $i,j\in\{1,2,\cdots,n\}$.
  \begin{align}
    \left|(H^{-1})_{i,j}\right| &\leq \frac{1}{\lambda_{\min}} \left(\frac{\lambda_{\max}-\lambda_{\min}}{\lambda_{\max}+\lambda_{\min}}\right)^{d_G(i,j)/\mathcal{B}_G(H)}      \label{eqn:spec}
  \end{align}
\end{thm}
The proof is given in Appendix \ref{thm:domdec-pf}. Theorem \ref{thm:domdec-gen} states that, when the graph-structured matrix $H$ is inverted, the magnitude of the $(i,j)$th component of its inverse decreases exponentially with respect to the distance between vertices $i$ and $j$ in $G$. Furthermore, the exponential decrease rate depends on the conditioning of $H$. This exponential decrease can be interpreted as a {\em diffusive} process that captures long-range interactions along vertices in the graph. {In Appendix \ref{apx:generalize}, we present a generalization of Theorem \ref{thm:domdec} that does not require $H$ to be PD.}

{
  \begin{rmk}
    In the remainder of the paper we focus on solving problems of the form \eqref{eqn:lin} (or equivalently, \eqref{eqn:elim}). Specifically, the problem of interest is $Hx=f$ with graph-structured $H$. Problem \eqref{eqn:orig-graph} is used to motivate our setting, but the problem of interest is not required to originate from \eqref{eqn:orig-graph}.
  \end{rmk}
}

\section{Decentralized Architecture with Overlap}\label{sec:dec}
{In this section we analyze the convergence of a decentralized algorithm for solving graph-structured optimization problems. We first focus on solving a graph-structured linear system \eqref{eqn:lin} with synchronous (Section \ref{sec:sync}) and asynchronous schemes (Section \ref{sec:async}), and then discuss its application to more general optimization settings (Section \ref{sec:opt}).} 

We consider the linear system \eqref{eqn:lin} with PD matrix $H\in\mathbb{R}^{n\times n}$ and associated  graph $G(V,E)$ with $V=\{1,2,\cdots,n\}$. To describe the decomposition scheme, we consider a partition of the entire domain $V$ into subdomains $\{V_k:k\in\mathcal{K}\}$, where $\mathcal{K}:=\{1,2,\cdots,K\}$ denotes the set of subdomains. Since we can relabel the nodes we can assume (without loss of generality) that:
\begin{align}
  V_k :=\left\{\sum_{k'=1}^{k-1}|V_{k'}|+1,\sum_{k'=1}^{k-1}|V_{k'}|+2,\cdots,\sum_{k'=1}^{k}|V_{k'}|\right\}, 
\end{align}
for $k\in\mathcal{K}$. For each subdomain $k \in \mathcal{K}$, we consider an overlapping region with neighboring subdomains of the form:
\begin{align}\label{eqn:ovlblk}
  V^\omega_k := \{i\in V \mid d_G(i,V_k) \leq \omega \}.
\end{align}
Note that $\exists\, k'\neq k$ such that $V^\omega_k \cap V^\omega_{k'}\neq\emptyset$ for any $k\in\mathcal{K}$ if $\omega\geq 1$ and the graph is connected. 

\subsection{Synchronous Coordination}\label{sec:sync}
We now describe a synchronous scheme to solve \eqref{eqn:lin}. The scheme can be interpreted as a {\em block-Jacobi} algorithm. For convenience, we define the projections of $H$ and $f$:
\begin{subequations}\label{eqn:blkmat}
  \begin{align}
    &H^\omega_{k} := \{(H)_{i,j}\}_{i\in V^\omega_k, j\in V^\omega_k}= (T^\omega_k)^T H T^\omega_k \\
    &H^\omega_{-k} := \{(H)_{i,j}\}_{i \in V^\omega_k, j \in V\setminus V^\omega_k}= (T^\omega_k)^T H T^\omega_{-k}\\
    &f^\omega_{k} := \{f\}_{V^\omega_k}= (T^\omega_k)^T f ,
  \end{align}
\end{subequations}
where $T_k :=\{e_i\}_{:,i\in V_k}$, $T^\omega_k := \{e_i\}_{:,i\in V^\omega_k}$, $T^\omega_{-k} := \{e_i\}_{:,i\in V\setminus V^\omega_k}$, and $e_i$ are standard unit vectors. If $V^\omega_k=V$, so $V\setminus V^\omega_k$ is empty, we consider $H^\omega_{-k}$ and $T^\omega_{-k}$ as empty matrices. To avoid confusion, we use the original indices after the projection (e.g., $(H^\omega_{k})_{i,j} = (H)_{i,j}$ for $i,j\in V^\omega_k$). Note that $H^\omega_{k}$ is also PD. The linear system \eqref{eqn:lin} can be solved in a decentralized manner as follows. First consider a system where, at iterate $t$, only the subset $\{x\}_{V^\omega_k}$ of variables is updated and the rest of the variables $\{x\}_{V\setminus V^\omega_k}$ are fixed to their current values. The corresponding subsystem has the form
\begin{align}\label{eqn:jacobi-sub}
  (H^\omega_{k})\{x\}_{V^\omega_k} =  \left(-H^\omega_{-k}\{x^{(t)}\}_{V\setminus V^\omega_k} + f_k^\omega\right),
\end{align}
where the solution at iterate $t$ is denoted by $x^{(t)}\in\mathbb{R}^n$.
After solving \eqref{eqn:jacobi-sub}, the solution is restricted to $V_k$ to obtain the next iterate. This restriction can be expressed as
\begin{align}\label{eqn:jacobi-indwise}
  x_k^{(t+1)}
  &= \left\{(H^\omega_{k})^{-1} (-H^\omega_{-k}\{x^{(t)}\}_{V\setminus V^\omega_k} + f_k^\omega)\right\}_{V_k}
\end{align}
for each subsystem $k\in\mathcal{K}$. Note that \eqref{eqn:jacobi-indwise} for each subsystem $k$ can be solved in parallel. By solving \eqref{eqn:jacobi-indwise}, one constructs the entire iterate vector $x^{(t+1)}=(x_1^{(t+1)},\cdots,x_K^{(t+1)})$. 

The iteration scheme can be expressed in the following matrix form:
\begin{align}\label{eqn:jacobik}
  x_k^{(t+1)} = S_k^\omega x^{(t)} + U^\omega_k f ,
\end{align}
where
\begin{subequations}\label{eqn:SUk}
  \begin{align}
    S^\omega_k &:= -(T_k)^T T^\omega_k (H^\omega_{k})^{-1} H^\omega_{-k}(T^\omega_{-k})^T\label{eqn:Sk}\\
    U^\omega_k &:=   (T_k)^T T^\omega_k (H^\omega_{k})^{-1} (T^\omega_k)^T \label{eqn:Uk}.
  \end{align}  
\end{subequations}
{The definition of $S^\omega_k$ \eqref{eqn:Sk} can be interpreted as (i) mapping the full-dimension vector $x$ onto $V\setminus V^\omega_k$ (to form the boundary information for $V^\omega_k$), (ii) multiplying the off-diagonal block $H^\omega_{-k}$ (this can be interpreted as the interaction with neighboring subdomains), (iii) solving the subproblem, (iv) mapping the solution on $V^\omega_k$ to the full space, and (v) restricting the solution to $V_k$. Similarly, the definition of $U^\omega_k$ \eqref{eqn:Sk} can be interpreted as (i) mapping the full-dimension vector to $V^\omega_k$, (ii) solving the subproblem, (iii) mapping the solution on $V^\omega_k$ to the full space, and (iv) restricting the solution to $V_k$.}

The overall iteration can be expressed as the linear system
\begin{align}\label{eqn:jacobi}
  x^{({t}+1)} = S^\omega x^{(t)} + U^\omega f ,
\end{align}
where
\begin{align}\label{eqn:SU}
  S^\omega :=[(S^\omega_1)^T  \cdots  (S^\omega_K)^T ]^T, \ U^\omega := [ (U^\omega_1)^T \cdots (U^\omega_K)^T ]^T.
\end{align}
We call \eqref{eqn:jacobi} a synchronous coordination scheme for \eqref{eqn:lin}.

We now derive conditions under which the synchronous scheme converges to the solution of \eqref{eqn:lin}.
\begin{lem}\label{lem:conv}
  Consider \eqref{eqn:lin} with PD $H\in\mathbb{R}^{n\times n}$ and associated graph $G(V,E)$ with $V=\{1,2,\cdots,n\}$. Consider also the synchronous scheme  \eqref{eqn:jacobi} with associated partitions $\{V_k\}_{k\in\mathcal{K}}$ and overlapping regions $ \{V^\omega_k\}_{k\in\mathcal{K}}$ with $\mathcal{K}:=\{1,2,\cdots,K\}$ and $\omega\geq 0$. The following statements are equivalent.
  \begin{enumerate}[label=(\alph*)]
  \item\label{lem:conv-a} $I-S^\omega$ is nonsingular, and the sequence generated by \eqref{eqn:jacobi} converges to the solution of \eqref{eqn:lin} as ${t}\rightarrow\infty$ (for any $x^{(0)}$).
  \item\label{lem:conv-b} $\rho(S^\omega)<1$
  \end{enumerate}
\end{lem}
The proof is given in Appendix \ref{lem:conv-pf}. 

The spectral radius $\rho(S^\omega)$ satisfies \cite[Theorem 5.6.9]{horn1990matrix}:
\begin{align}
  \rho(S^\omega)\leq \Vert S^\omega \Vert_\infty.
\end{align}
Thus, we need only to find conditions under which $\Vert S^\omega\Vert_\infty<1$ in order to guarantee convergence. The following lemma establishes such conditions.

\begin{lem}\label{lem:decen}
  Consider \eqref{eqn:lin} with PD $H\in\mathbb{R}^{n\times n}$ and the associated graph $G(V,E)$ with $V=\{1,2,\cdots,n\}$. Consider also the synchronous scheme  \eqref{eqn:jacobi} for a partition $\{V_k\}_{k\in\mathcal{K}}$ of $V$ and overlapping regions $\{V^\omega_k\}_{k\in\mathcal{K}}$ with $\mathcal{K}:=\{1,2,\cdots,K\}$ and $\omega\geq 0$. We have that
  \begin{align}\label{eqn:Sombound}
    \Vert S^\omega\Vert_\infty \leq \max_{k\in\mathcal{K}} \frac{R_k}{\lambda^k_{\min}} \left(\frac{\lambda^k_{\max}-\lambda^k_{\min}}{\lambda^k_{\max}+\lambda^k_{\min}}\right)^{(\omega+1)/\mathcal{B}_G(H^\omega_{k})-1} 
  \end{align}
  where
  $\lambda^k_{\max}$ and $\lambda^k_{\min}$ are the largest and the smallest eigenvalues of $H^\omega_{k}$, respectively; and $R_k:=\sum_{i\in V^\omega_k,j\in V\setminus V^\omega_{k}} \left|(H)_{i,j}\right|$.
\end{lem}
The proof is given in Appendix \ref{prop:decen-pf}.

{We now analyze the {\em convergence rate} of the coordination scheme. We note that, from \eqref{eqn:jacobi}, we have that 
  \begin{align}\label{eqn:pre-main-1}
    x^{(t)}-x^* =(S^\omega)^t \left(x^{(0)} - x^*\right).
  \end{align}  
  By taking $\Vert\cdot\Vert_\infty$ on the both side of \eqref{eqn:pre-main-1}, we have:
  \begin{align}\label{eqn:pre-main-2}
    \left\Vert x^{(t)}-x^*\right\Vert_\infty
    &\leq \Vert S^\omega\Vert_\infty^t \left\Vert x^{(0)}-x^*\right\Vert_\infty.
  \end{align}
Thus,   if $\Vert S^\omega\Vert_\infty <1$, the convergence is linear. 
  The following theorem is the result of Lemma \ref{lem:decen} and \eqref{eqn:pre-main-2}.
  \begin{thm}\label{thm:main}
    Consider \eqref{eqn:lin} with PD $H\in\mathbb{R}^{n\times n}$ and the associated graph $G(V,E)$ with $V=\{1,2,\cdots,n\}$. Consider also the synchronous scheme  \eqref{eqn:jacobi} for a partition $\{V_k\}_{k\in\mathcal{K}}$ of $V$ and overlapping regions $\{V^\omega_k\}_{k\in\mathcal{K}}$ with $\mathcal{K}:=\{1,2,\cdots,K\}$ and $\omega\geq 0$. The sequence generated by \eqref{eqn:jacobi} satisfies:
    \begin{align}\label{eqn:main}
      \left\Vert x^{(t)}-x^*\right\Vert_\infty \leq \alpha^t \left\Vert x^{(0)} - x^*\right\Vert_\infty
    \end{align}
    where the convergence rate $\alpha$ is defined by:
    \begin{align}\label{eqn:alpha}
      \alpha:= \max_{k\in\mathcal{K}} \frac{R_k}{\lambda^k_{\min}} \left(\frac{\lambda^k_{\max}-\lambda^k_{\min}}{\lambda^k_{\max}+\lambda^k_{\min}}\right)^{(\omega+1)/\mathcal{B}_G(H^\omega_{k})-1}.
    \end{align}
\end{thm}}

The convergence rate is bounded by $\Vert S^\omega\Vert_\infty$ and decreases (improves) exponentially with the size of the overlap $\omega$. 
\begin{rmk}\label{rmk:further}
The bound of convergence rate $\alpha$ can be further simplified to
  \begin{align}
    \alpha \leq \frac{R}{\lambda_{\min}}\left(\frac{\lambda_{\max}-\lambda_{\min}}{\lambda_{\max}+\lambda_{\min}}\right)^{(\omega+1)/\mathcal{B}_{G}(H)-1} ,
  \end{align}
  where $R:=\max_{k\in\mathcal{K}} R_k$, $\lambda_{\max}$ and $\lambda_{\min}$ are the largest and smallest eigenvalues of $H$, respectively.
\end{rmk}
{\begin{rmk}\label{rmk:pde}
    A similar form of the bound of convergence rate \eqref{eqn:alpha} is also observed in the PDE literature. Specifically, the overlapping Schwarz algorithm, when applied to elliptic PDEs, has a similar form of convergence rate; the rate decreases exponentially with the size of overlap, and the decrease rate is related to the eigenvalues of the problem \cite{dolean2015introduction}.
\end{rmk}}
\begin{rmk}\label{rmk:offdiagonal}
  Condition \eqref{eqn:Sombound} shows that having a small $R_k$ is also important in order to achieve fast convergence. Intuitively, $R_k$ can be interpreted as the magnitude of the coupling between $V_k^\omega$ and $V\setminus V_k^\omega$. Thus, when partitioning a graph, such information should be taken into account (the weaker the coupling, the faster the convergence).
\end{rmk}
\begin{rmk}\label{rmk:tradeoff}
While the convergence rate improves with larger overlap $\omega$, the computational cost for the subdomain problems increases with $\omega$ (because the size of each overlapping block $V_k^\omega$ increases with $\omega$). Also, information sharing requirements increase with $\omega$. This situation highlights the inherent trade-off between performance and data privacy. The size of the overlap $\omega$ and the partitioning strategy are key design parameters that can help trade off these aspects and cover an entire spectrum of architectures that have fully centralized and fully decentralized architectures as extreme points. Specifically, a fully centralized scheme is obtained when the size of the overlap is the dimension of the entire system, whereas a fully decentralized scheme is that in which the overlap is zero.
\end{rmk}

\subsection{Asynchronous Coordination}\label{sec:async}
We now analyze the convergence properties of the architecture under asynchronous coordination. We consider a set of iteration indices $\mathcal{T}:=\{0,1,2,\cdots\}$ at which one or multiple subdomains are updated. Also, we denote the set of times at which $x_k$ is updated by $\mathcal{T}_k$. We assume that when updating $x_k$, the information for the other subdomains may not be the most recent ones. That is, rather than using the information $x^{(t)}$, it uses a combination of the present and delayed information
$x^{k,(t)} := (x_1^{(\tau_{k,1}(t))},x_2^{(\tau_{k,2}(t))},\cdots,x_K^{(\tau_{k,K}(t))})$  to compute $x_k^{(t+1)}$. Here, $\tau_{k,k'}(t)$ is the (delayed) time index of the information from subdomain $k'$ used to compute $x_k^{(t+1)}$. The asynchronous iteration takes the  form
\begin{align}\label{eqn:async}
  x_k^{({t}+1)} =
  \begin{cases}
    S^\omega_k
    x^{k,(t)}+ U^\omega_k f
    & t\in\mathcal{T}_k\\
    x_k^{(t)} &\text{otherwise} 
  \end{cases}
\end{align}
for $k\in\mathcal{K}$ and $t\in\mathcal{T}$.
The asynchronous scheme becomes the synchronous counterpart whenever $\tau_{k,k'}(t)=t$ for all $k,k'\in\mathcal{K}$ and $t\in\mathcal{T}$ and $\mathcal{T}_1=\cdots=\mathcal{T}_k=\mathcal{T}$. In our analysis we make the following assumptions \cite[Assumption 1.1, pg 430]{bertsekas1989parallel}.

\begin{ass}[Total Asynchronism]\label{ass:totalasync}
  For $\mathcal{T}:=\{0,1,2,\cdots\}$, $\mathcal{T}_k\subseteq \mathcal{T}$ for $k\in\mathcal{K}$, and $0\leq \tau_{k,k'}(t)\leq t$ for $k,k'\in\mathcal{K},t\in\mathcal{T}_k$:
  \begin{enumerate}[label=(\alph*)]
  \item The sets $\mathcal{T}_k \subseteq \mathcal{T}$ are infinite for $k\in\mathcal{K}$.\label{ass:totalasync-b}
  \item $\lim_{t\rightarrow\infty} \tau_{k,k'}(t) = \infty$. \label{ass:totalasync-c}
  \end{enumerate}
\end{ass}
{In other words, the processor continues to iterate indefinitely and any old information is discarded in a finite number of iterations. Thus, Assumption \ref{ass:totalasync} holds unless one of the processors or the communication channels completely fails.}

The following provides a sufficient condition for the convergence of the asynchronous scheme \eqref{eqn:async}. 
\begin{thm}\label{thm:async}
  Consider \eqref{eqn:lin} with PD $H\in\mathbb{R}^{n\times n}$ and associated graph $G(V,E)$  with $V=\{1,2,\cdots,n\}$ as well as a partition $\{V_k\}_{k\in\mathcal{K}}$ of $V$ and overlapping regions $\{V^\omega_k\}_{k\in\mathcal{K}}$ with $\mathcal{K}:=\{1,2,\cdots,K\}$ and $\omega\geq 0$.  Consider also sets $\mathcal{T}$, $\mathcal{T}_k$, and $\tau_{k,k'}(t)$ for $k,k'\in\mathcal{K},t\in\mathcal{T}_k$ such that Assumption \ref{ass:totalasync} holds. The asynchronous scheme \eqref{eqn:async} converges to the solution of \eqref{eqn:lin} if $\Vert S^\omega \Vert_\infty<1$.
\end{thm}
The proof is given in Appendix \ref{prop:async-pf}. The proof is based on \cite[Proposition 2.1, pg 431]{bertsekas1989parallel}. Lemma \ref{lem:conv} and Theorem \ref{thm:async} establishes that convergence is guaranteed if $\omega$ is sufficiently large.

\subsection{Operating in the Optimization Space}\label{sec:opt}
The proposed scheme can be implemented by operating in the space of the optimization problem. To see this, consider
\begin{align}\label{eqn:opt}
  \min_{x\in\mathbb{X}} \quad \varphi(x),
\end{align}
where $\mathbb{X}\subseteq\mathbb{R}^n$ is closed and $\varphi:\mathbb{R}^n\rightarrow \mathbb{R}$. We also consider a graph $G(V,E)$ with $V=\{1,2,\cdots,n\}$, partitions $\{V_k\}_{k\in\mathcal{K}}$ with $\mathcal{K}=\{1,2,\cdots,K\}$, and its overlapping blocks $\{V^\omega_k\}_{k\in\mathcal{K}}$ with $\omega\geq 0$. In the synchronous scheme for the solution of \eqref{eqn:opt}, the subproblem can be written as
\begin{subequations}
  \begin{align}\label{eqn:opt-sub}
    \min_{ x\in\mathbb{X}}\quad &\varphi(x)\\
    \text{s.t.}\quad &\{x\}_{V\setminus V^\omega_k}=\{x^{(t)}\}_{V\setminus V^\omega_k}\label{eqn:opt-sub-con},
  \end{align}
\end{subequations}
where $t$ is the iteration counter. This is a counterpart of \eqref{eqn:jacobi-sub} in the optimization space. The solution is then restricted to the non-overlapping block $V_k$. The overall procedure can be written as follows:
\begin{align}\label{eqn:opt-jacobi}
  x_k^{(t+1)} = \left\{\mathop{\text{argmin}}_{x\in\mathbb{X}} \varphi(x)\; \text{s.t.}\; \eqref{eqn:opt-sub-con}\right\}_{V_k}.
\end{align}
This optimization subproblem corresponds to \eqref{eqn:jacobik}. The representation as optimization problems reveals that the proposed paradigm can, in principle, be applied to general problems.


\section{Implementation Details}\label{sec:impl}
In this section, we discuss the implementation of the algorithms presented in Section \ref{sec:dec} in a distributed computing environment. We first explain the implementation of the asynchronous scheme, and later briefly discuss the implementation of the synchronous scheme as a special case. Here we use MPI \cite{gropp1996high}, a portable message-passing standard. MPI enables efficient communication between parallel computing processors and can be used on a wide range of modern computer architectures. We assume that a processor and a memory are available for each subdomain $k\in\mathcal{K}$. These are referred to as process $k$ and memory $k$, respectively. Process $k$ performs computation and communication required for updating the solution and monitoring convergence of the algorithm. Communication between the processes is performed by using {\em one-sided communication} with remote memory access (RMA) operations \cite{gropp1996high}. In one-sided communication, a part of the memory in each process is declared as public memory region and exposed to RMA operations called from other processes. Accordingly, information that needs to be available to other processes is stored in the public region, while other information is stored in the private memory region.

The information that memory $k$ needs to contain is: 
\begin{itemize}
\item The private region of memory $k$ stores solution, $x^k=(x_1^k,x_2^k,\cdots,x_K^k)$,  errors $\epsilon^k=(\epsilon^k_1,\epsilon^k_2,\cdots,\epsilon^k_K)$,  matrices $S^\omega_k$, $H^\omega_{\pm k} $, and vectors  $U^\omega_k f$ and $f^\omega_k$
\item The public region of memory $k$ stores a copy of the local solution $\tilde{x}_k$ and local error $\tilde{\epsilon}_k$.
\end{itemize}
Here, $H^\omega_{\pm k}:=\{(H)_{i,:}\}_{i\in V^\omega_k}$. We use the superscript to highlight in which distributed memory the variable is stored, and we use the notation $\tilde{\cdot}$ to specify that the variable is stored in the public memory. Process $k$ performs the computations:
\begin{itemize}
\item Update the local solution $x_k^k \leftarrow S^\omega_k x^k + U^\omega_k f$.
\item Update the local error $\epsilon_k\leftarrow \Vert H^\omega_{\pm k} x^k - f^\omega_k\Vert_\infty $.
\end{itemize}
After the computation, the information on the states and the errors is exchanged across the parallel processes, and the procedure is repeated. The following remark discusses an alternative implementation.
\begin{rmk}\label{eqn:invsps}
The computations for the synchronous \eqref{eqn:jacobi} or asynchronous scheme \eqref{eqn:async} can be performed in two ways. One can explicitly compute $S^\omega_k$, which includes explicitly computing $(H^\omega_{k})^{-1}$, in advance and simply perform matrix multiplications in each iteration. Alternatively, one can solve sparse systems \eqref{eqn:jacobi-sub} instead of inverting $H^\omega_{k}$. The first strategy is more practical when the problem needs to be repeatedly solved with a fixed $H$ and varying $f$ (e.g., control and estimation);  the latter is more effective when $H_{k}^\omega$ changes or is difficult to be inverted.
\end{rmk}

\begin{figure}[h!]
  \centering
  \includegraphics[width=.3\textwidth]{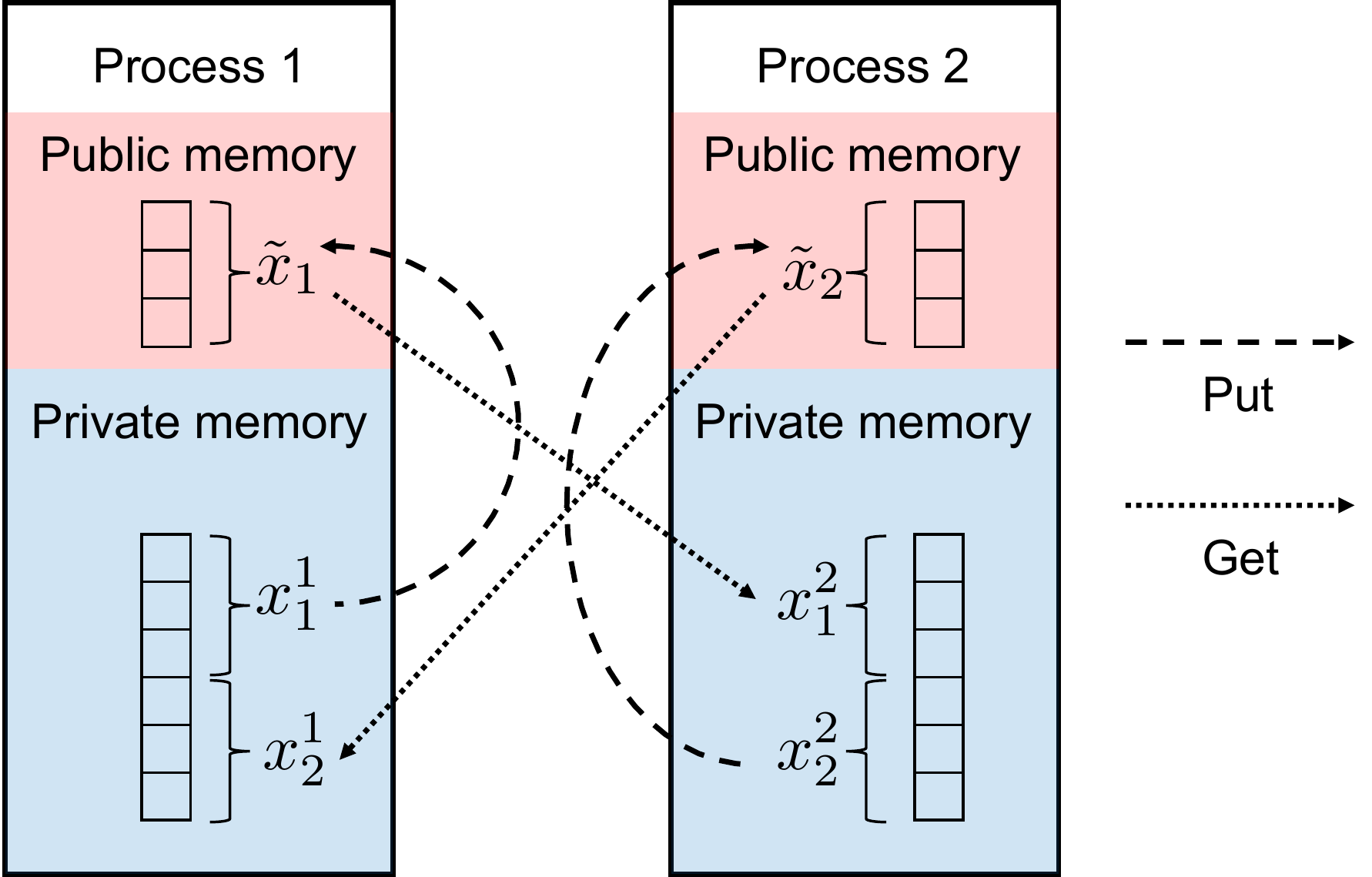}
  \caption{Schematic of parallel implementation with RMA operations.}\label{fig:impl}
\end{figure}

The communication between processes is performed as follows. Each process puts the local information (solution and error) on the public window of its local memory and  gets the local information (solution and error) of other processes by remotely accessing the public region of the other processes memory (see Fig. \ref{fig:impl}). In particular, such communication is performed by RMA functions {\tt MPI\_Win\_put} and {\tt MPI\_Win\_get}. The communication procedure can be described as follows.
\begin{itemize}
\item Process $k$ {\em puts} the local solution $x^k_k$ and the local error $\epsilon^k_k$ to the public window (i.e., $\tilde{x}_k\leftarrow x^k_k$ and $\tilde{\epsilon}_k\leftarrow \epsilon^k_k$).
\item Process $k$ {\em gets} the local solution $\tilde{x}_{k'}$ and the local error $\tilde{\epsilon}_{k'}$ of process $k'$ from the public window of memory $k'$ for each $k'\neq k$ (i.e., $x^{k}_{k'} \leftarrow \tilde{x}_{k'}$ and $\epsilon^k_{k'}\leftarrow \tilde{\epsilon}_{k'}$).
\end{itemize}
Simultaneous access to a single memory can cause an issue because information can be updated while being read by other processes. In order to prevent such an issue, the public region of the memory is locked while the memory is being used. While a public memory is being read (by RMA get operations), only simultaneous readings are allowed. While a public memory is being written (by RMA put operations), simultaneous access of any kind (either reading or writing) is not allowed. Such features can be implemented by using RMA synchronization operations. In particular, {\tt MPI\_Win\_lock} and {\tt MPI\_Win\_unlock} are used.

The overall implementation is summarized in Algorithm \ref{alg:decen}. The only difference between the synchronous and asynchronous mode is that in synchronous mode the processes are synchronized before and after putting the solution $x_k$ to the public memory of process $k$. This approach enforces the iteration of each process to proceed at the same speed.

\begin{algorithm}
  \caption{Implementation of Decentralized Architecture}\label{alg:decen}
  \begin{algorithmic}
    \For {$k=1,\cdots,K$ (In parallel)}
    \State Initialize $x_k\leftarrow 0$, $\epsilon \leftarrow \infty$
    \While {$\epsilon>\epsilon^{\text{tol}}$}
    \State {\bf *Synchronize} processes $k=1,\cdots,K$
    \State \textbf{Lock (exclusive)} memory $k$
    \State \textbf{Put} $\tilde{x}_k\leftarrow x_k$ to memory $k$
    \State \textbf{Put} $\tilde{\epsilon}_k\leftarrow \epsilon_k$ to memory $k$
    \State \textbf{Unlock (exclusive)} memory $k$
    \State {\bf *Synchronize} processes $k=1,\cdots,K$
    \For {$k'=1,\dots,K$ and $k'\neq k$}
    \State \textbf{Lock (shared)} memory $k$
    \State \textbf{Get} $x^k_{k'}\leftarrow\tilde{x}_{k'}$ from memory $k'$
    \State \textbf{Get} $\epsilon^k_{k'}\leftarrow\tilde{\epsilon}_{k'}$ from memory $k'$
    \State \textbf{Unlock (shared)} memory $k$
    \EndFor
    \State Update $\epsilon_k\leftarrow \Vert H^\omega_{\pm k} x^k - f^\omega_k \Vert_\infty$
    \State Update $x_k^k \leftarrow S^\omega_k x^k + U^\omega_k f$
    \State $\epsilon\leftarrow \max_{k'\in\mathcal{K}}\{ \epsilon^k_{k'} \}$ 
    \EndWhile 
    \EndFor
  \end{algorithmic}
  *Performed only in synchronous mode.
\end{algorithm}

\section{Case Study}\label{sec:cstudy}
We demonstrate the proposed framework by solving a state estimation problem for a 9241-bus power network. We consider a power network system on a graph $G(V,E)$. We assume that the network is primarily inductive, the voltage amplitudes are fixed to one, and the voltage angle differences between the neighboring nodes are small enough to apply a DC approximation. The power flow $P_{ij}$ on edge $\{i,j\}\in E$ can be expressed by $P_{ij} = y_{ij} (\delta_i - \delta_j)$ (assume that a direction is assigned to each edge). We assume that the power flow is measured and the measurement is performed based on a statistical model $P^m_{ij} = P_{ij}+\eta_{P_{ij}}$, where $\eta_{P_{ij}}$ is a random variable whose distribution is $\eta_{P_{ij}}\sim N(0,\sigma^2_{P_{ij}})$. By incorporating the prior on $\delta_i \sim N(\delta_i^m,\sigma_{\delta_i}^2)$ for $i\in V$, one can derive the following maximum a posteriori problem.
\begin{subequations}\label{eqn:est}
  \begin{align}
    \begin{split}\min_{\delta,P} \ &\sum_{i \in V} \left(\frac{\delta_i-\delta_i^m}{\sigma_{\delta_i}}\right)^2+ \sum_{\{i,j\}\in E} \left(\frac{P_{ij}-P_{ij}^m}{\sigma_{P_{ij}}}\right)^2
    \end{split}
    \\
    \text{s.t.}\ & P_{ij} =  y_{ij}(\delta_i - \delta_j), \quad \{i,j\}\in E
  \end{align}
\end{subequations}
In our estimation setting, we assume that only a subset of flows can be measured and the rest need to be inferred from data.
Accordingly, we assume $\sigma_{P_{ij}}=y_{ij}$ for the measured flows (about half of the edges are randomly selected) and assume much weaker prior for the rest of the edges by $\sigma_{P_{ij}}=\sqrt{10} y_{ij}$.
  The prior weight $c\in \mathbb{R}_{>0}$ on the unmeasured voltage angles is assumed to be uniform, that is, $\sigma_{\delta_i}=\frac{1}{c}, \forall i \in V$.  The estimation problem can be written in vector form as
\begin{subequations}\label{eqn:est-vec}
  \begin{align}
    \min_{\delta,P} \ & (\delta-\delta^m)^T \Sigma_\delta (\delta-\delta^m) + ({P} -{P}^m)^T \Sigma_{{P}} ({P}-{P}^m)\\
    \text{s.t.}\ & P = Y \delta ,
  \end{align}
\end{subequations}
where $Y\in\mathbb{R}^{|E|\times|V|}$ and $\Sigma_\delta,\Sigma_P \in\mathbb{R}^{|V|\times|V|}$. This problem can be reduced to
\begin{align}
  \min_{\delta} \ & \delta^T \left(\Sigma_\delta + Y^T \Sigma_{P}Y  \right)\delta - 2(Y^T \Sigma_P P^m +  \Sigma_{\delta}\delta^m )^T  \delta.
\end{align}
Here, we note that the assignments
$H:=\Sigma_\delta + Y^T \Sigma_{P}Y$, $f:= Y^T\Sigma_P P^m +  \Sigma_{\delta}\delta^m $, and $x:=\delta$ fit this problem into the paradigm of \eqref{eqn:lin}.
We have that $\mathcal{B}_G(\Sigma_\delta)=1$, $\mathcal{B}_G(Y^T \Sigma Y)\leq1$, and thus $\mathcal{B}_G(H) \leq 1$.

We used data from the Pegase project \cite{fliscounakis2013contingency} to derive the power system model. We applied graph partitioning based on a multilevel $k$-way partitioning method using {\tt METIS} \cite{karypis1995metis} to identify the partition $\{V_k\}_{k\in\mathcal{K}}$. Our implementation of Algorithm \ref{alg:decen} uses the popular {\tt MPICH} MPI library and the basic linear algebra subprograms (BLAS) package for matrix computation. The program was run on a multicore parallel computing server (four nodes and one CPU core Intel Xeon Processor E5-2695v4 per node). 

\begin{figure}[!t]
  \centering
  \includegraphics[width=.35\textwidth]{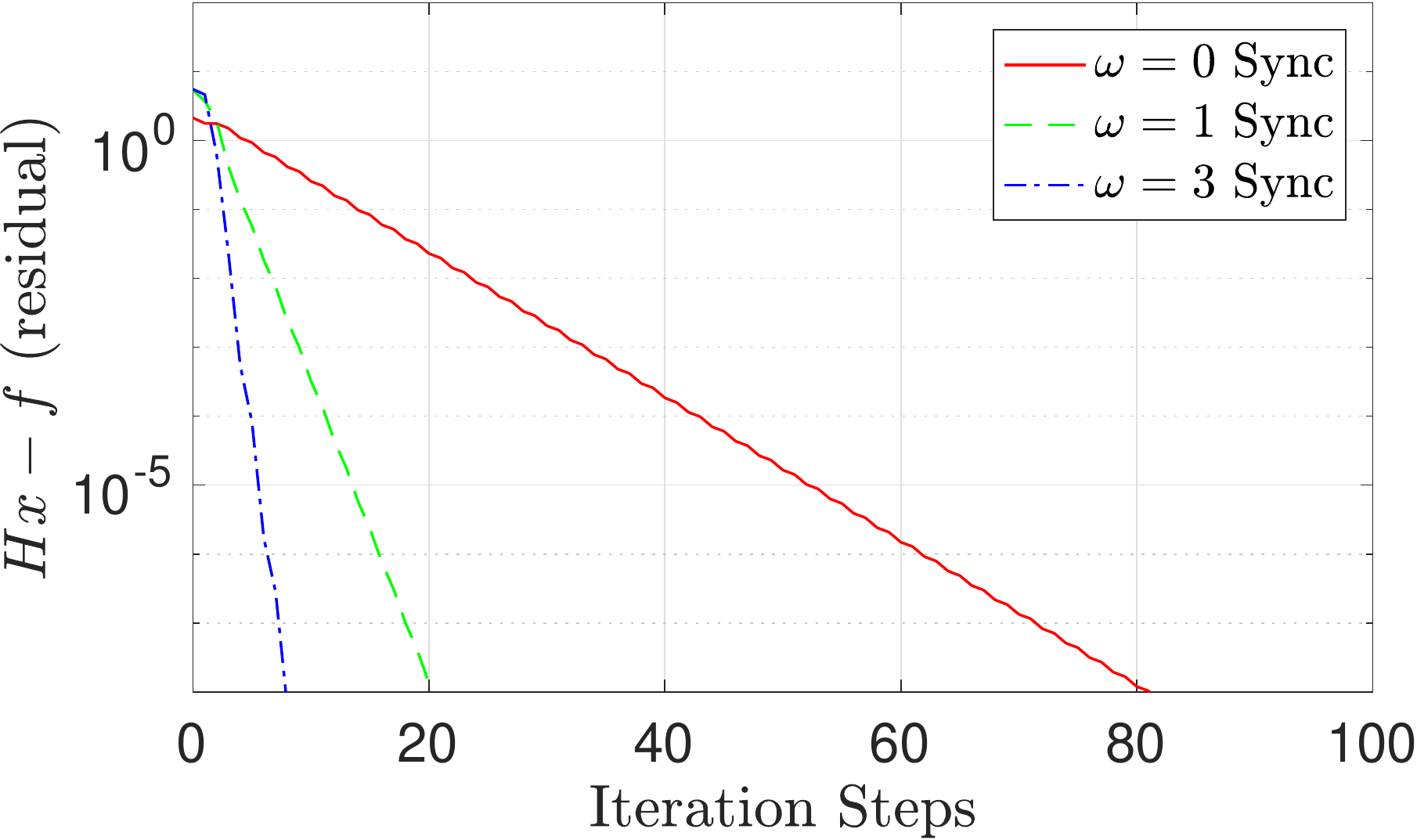}\\
  \includegraphics[width=.35\textwidth]{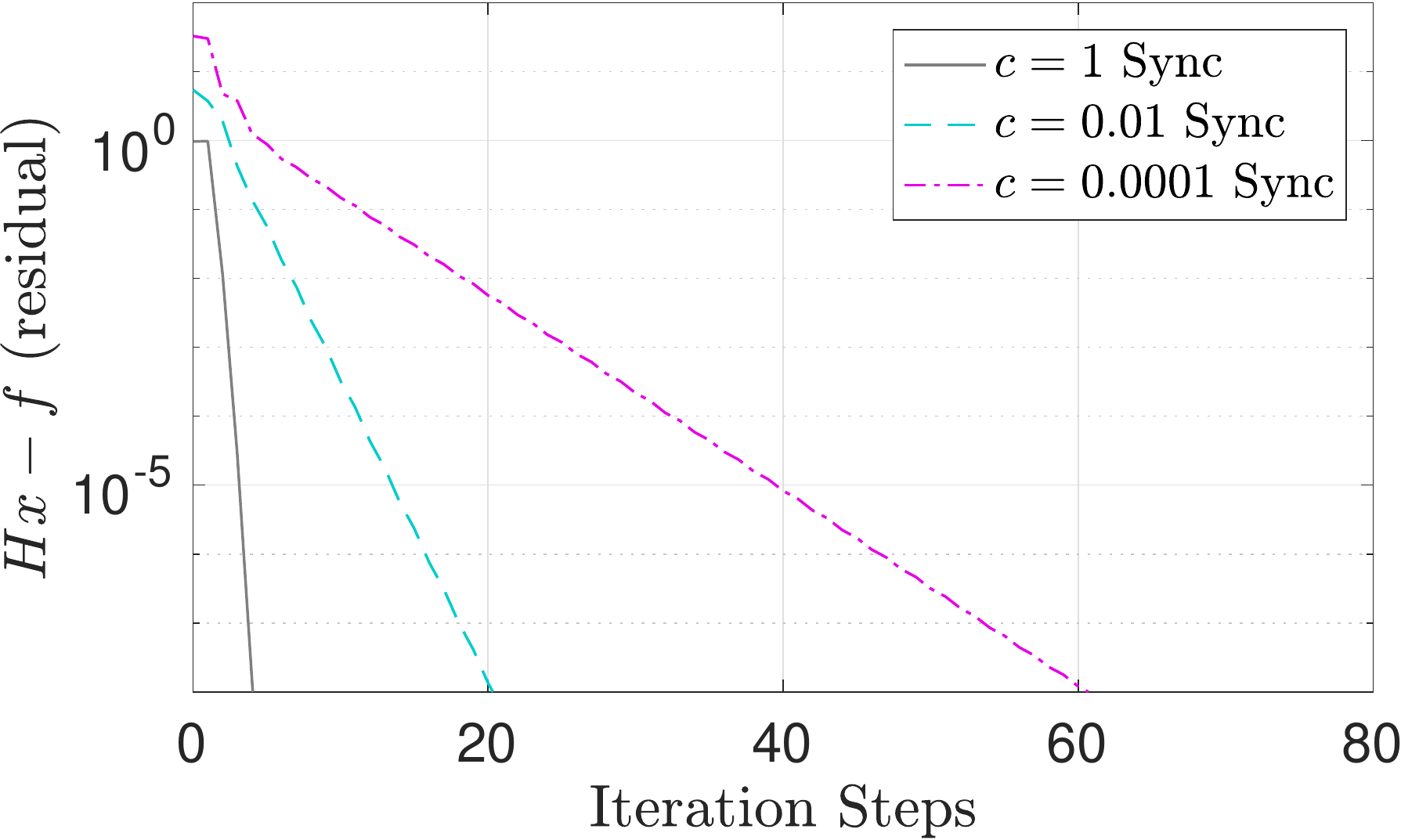}
  \caption{(Top) Residual over iteration steps ($c=0.1$, partition structure 1). (Bottom) Residual over iteration steps ($\omega=1$, partition structure 1).}\label{fig:res-iter}
\end{figure}

In Fig. \ref{fig:res-iter} (top), the residual $Hx^{(t)}-f$ for the synchronous mode is plotted with different sizes of the overlap. The results confirm that the solution converges linearly to $x^*=H^{-1}f$ and that it converges faster as the size of the overlap increases (see Lemma \ref{lem:decen}). In Fig. \ref{fig:res-iter} (bottom), we present the evolution of the residual for different values of the regularization coefficient $c$. If $c$ is large, $\lambda_{\min}$ increases and $\lambda_{\min}+\lambda_{\max}$ increases. We thus confirm (as shown in Lemma \ref{lem:decen}) that the convergence rate improves with larger $c$. 

\begin{figure}[!t]
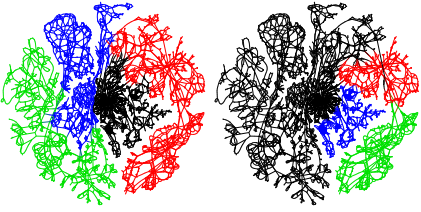

  \centering
  \includegraphics[width=.15\textwidth]{pegasus-part-4.pdf}
  \includegraphics[width=.15\textwidth]{pegasus-part-4-2.pdf}
  \caption{(Left) Network topology with partition structure 1. (Right) Network topology with partition structure 2.}
  \label{fig:network}
\end{figure}

\begin{figure}[!t]
  \centering
  \includegraphics[width=.35\textwidth]{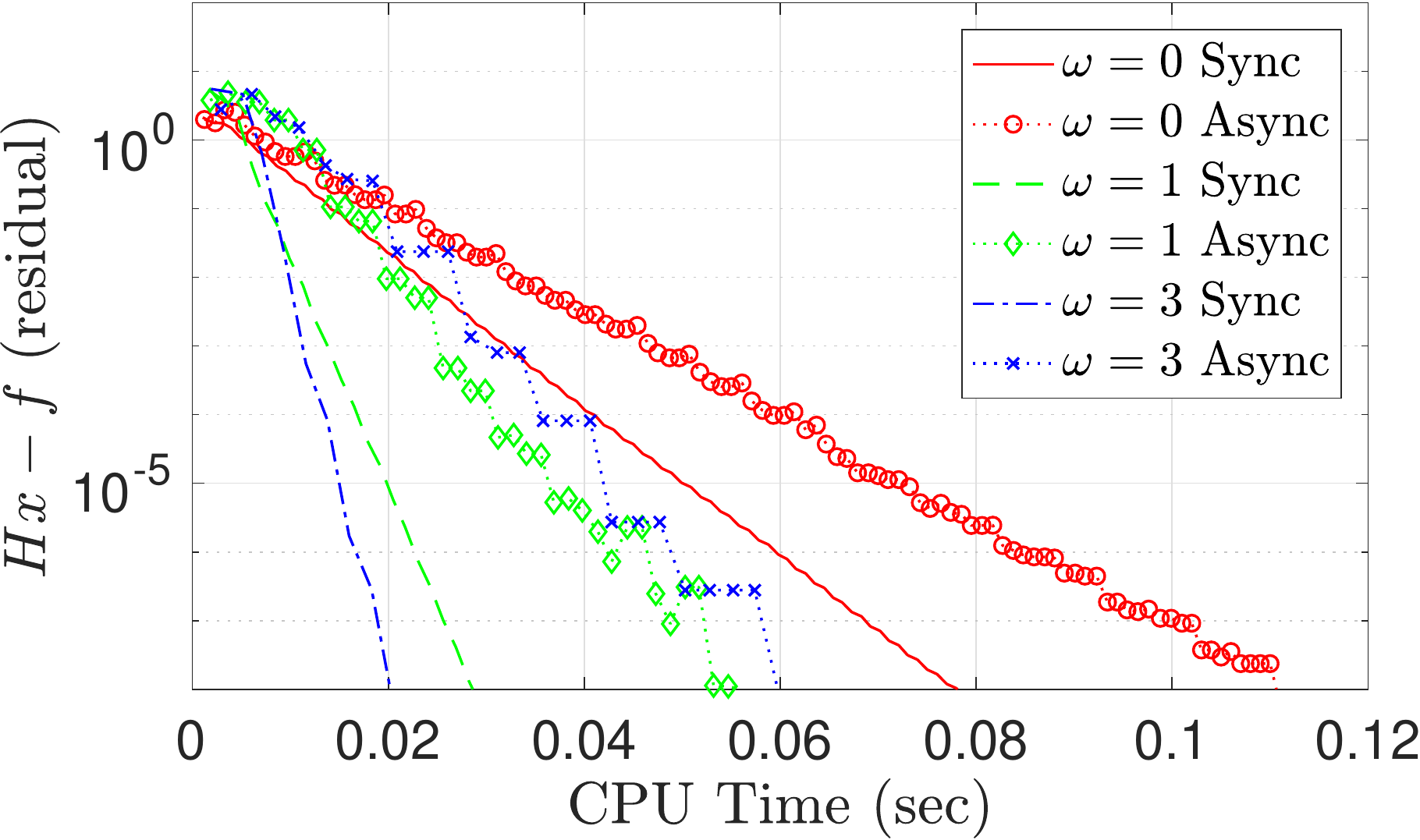}\\
  \includegraphics[width=.35\textwidth]{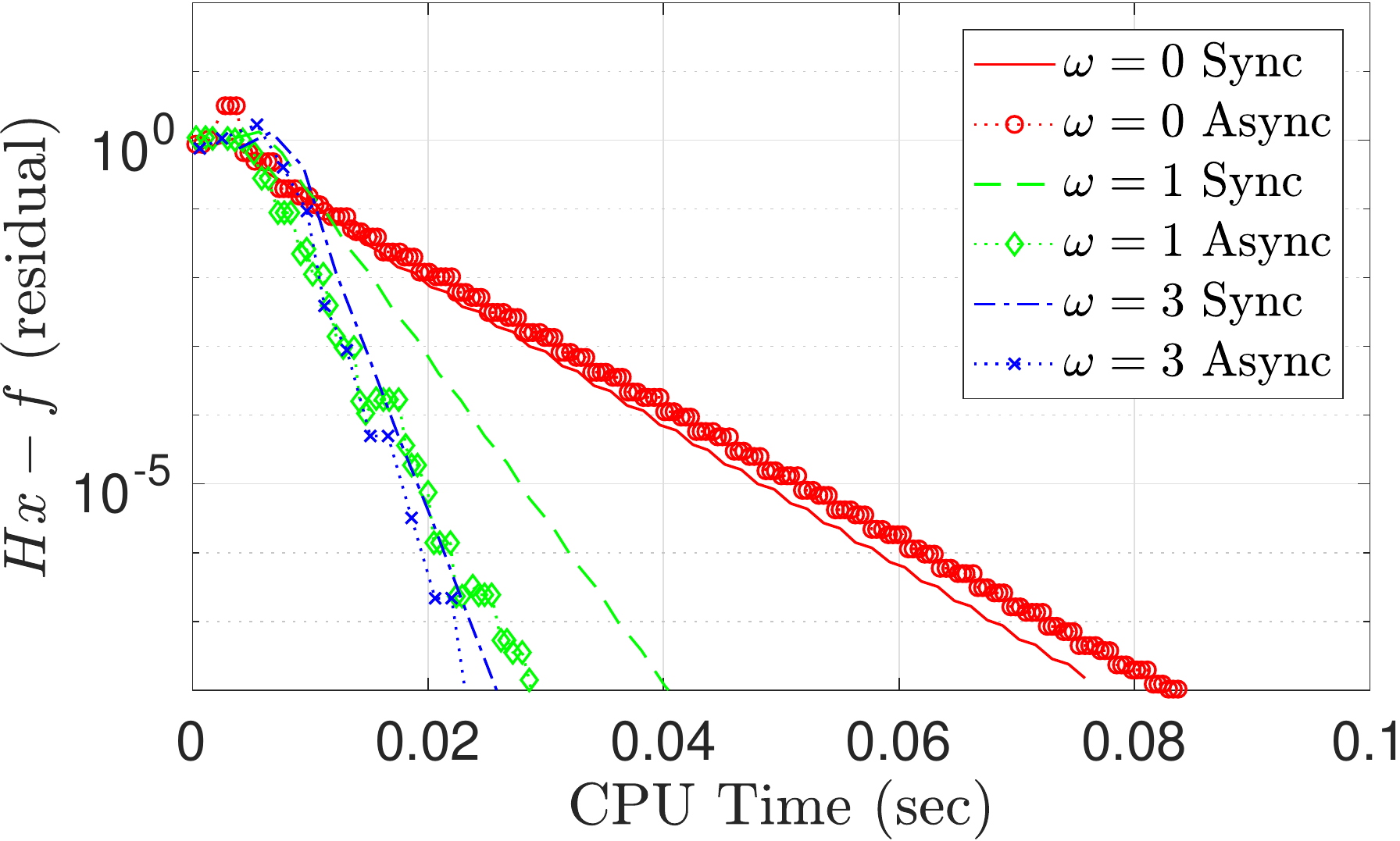}
  \caption{(Top) Residual over CPU time ($c=0.1$, partition structure 1). (Bottom) Residual over CPU time ($c=0.1$, partition structure 2). }\label{fig:res-cpu}
\end{figure}

In Fig. \ref{fig:res-cpu}, the residual for the synchronous and asynchronous implementations are contrasted. In Fig. \ref{fig:res-cpu} (top), we see that the synchronous scheme is faster than the asynchronous mode. The reason is that the synchronous scheme guarantees that each subproblem is solved with up-to-date information, whereas the asynchronous scheme does not. As can be observed in Fig. \ref{fig:res-cpu} (bottom), however, when a large imbalance exists in the computation and communication loads across the blocks (See Fig. \ref{fig:network}; {the dimension of the subdomains is uniform in partition structure 1 but non-uniform in partition structure 2}), the asynchronous scheme becomes faster than the synchronous scheme. The reason is that, in asynchronous mode,  processes with small blocks can iterate multiple times while the processes with large blocks perform a single iteration step (more effective use of computing resources is enabled). 

We also observe that the increase in the overlap reduces the number of iterations but not necessarily the solution time. The reason is that the increase in the overlap also increases the computation and communication cost (Remark \ref{rmk:tradeoff}). The overall computing cost increases with the size of the block (i.e., $|V^\omega_k|$), and the communication cost also tends to increase with the number of neighboring nodes (i.e., $|V_k^{\omega+1}\setminus V_k^\omega|$). We can also observe that the size of the overlapping blocks and the number of neighboring nodes increase with $\omega$ (Table \ref{tbl:stat}). Such effects ultimately are manifested in the overall CPU time per iteration: $0.0963$ sec/iter when $\omega=1$, $0.140$ sec/iter when $\omega=2$, and $0.256$ sec/iter when $\omega=3$. These results again illustrate the trade-offs in convergence and computational performance. 

\begin{table}[!t]
  \caption{Statistics for the subsystems (partition structure 1)}
  \label{tbl:stat}
  \centering
  \begin{tabular}{|c|c|c|c|c|c|}
    \hline
    k & 1 & 2 & 3 & 4 & Total\\ \hline
    $|V_k|$  & 2,324 & 2,366 & 2,278 & 2,273 & 9,241\\ \hline
    $|V_k^1|$  & 2,361 & 2,398 & 2,291 & 2,304 & 9,354\\ \hline
    $|V_k^2|$  & 2,452 & 2,469 & 2,322 & 2,376 & 9,619\\ \hline
    $|V_k^3|$  & 2,570 & 2,558 & 2,380 & 2,506 & 10,014\\ \hline
    $|V_k^4|$  & 2,744 & 2,680 & 2,485 & 2,727 & 10,636\\ \hline
    $|V_k^1\setminus V_k|$  & 37 & 32 & 13 & 31 & 113\\ \hline
    $|V_k^2\setminus V_k^1|$  & 91 & 71 & 31 & 72 & 265\\ \hline
    $|V_k^3\setminus V_k^2|$  & 118 & 89 & 58 & 130 & 395\\ \hline
    $|V_k^4\setminus V_k^3|$  & 174 & 122 & 105 & 221 & 622\\ \hline
  \end{tabular}
\end{table}

\begin{figure}[!t]
  \centering
  \includegraphics[width=.35\textwidth]{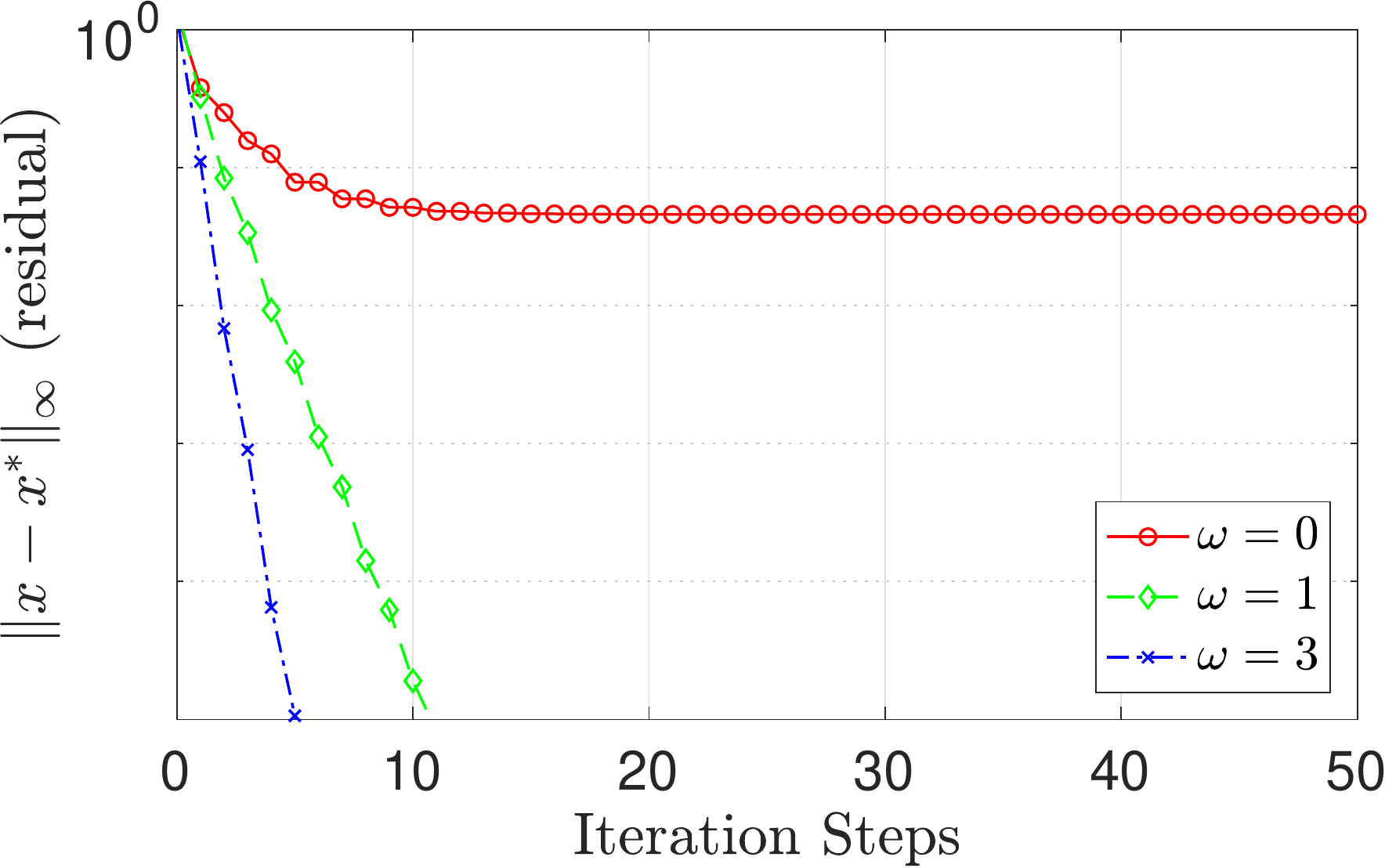}
  \caption{Error over iteration steps ($c=0.1$, partition structure 1) for the inequality-constrained optimization problem.}\label{fig:res-qp}
\end{figure}

We also performed computational experiments to analyze the performance of the algorithm in inequality-constrained QPs. Specifically, we incorporated the inequality constraints $-\pi/4\leq\delta_i-\delta_j\leq\pi/4\quad\text{for}\quad \{i,j\}\in E$.  We reformulated the problem using soft constraints and slack variables to make sure that the subproblem is always feasible. We applied the synchronous scheme \eqref{eqn:opt-jacobi} to solve the problem. The results are shown in Fig. \ref{fig:res-qp}. As can be seen, the algorithm converges when the overlap is sufficiently large (a fully decentralized approach with no overlap does not converge). This provides evidence that the proposed decentralized architecture with overlap can potentially be applied to solve large-scale constrained problems. We hypothesize that this is because the set of active constraints settles and the algorithm then reverts to a pure QP phase in a constrained subspace. We will investigate the convergence properties in this setting in future work. 

{We compared the performance of the proposed overlapping scheme with that of ADMM. To use ADMM, the original problem \eqref{eqn:elim} is rewritten in the lifted form:
\begin{subequations}\label{eqn:lifted}
  \begin{align}
    \min_{x,z} \quad& \sum_{k\in\mathcal{K}} \frac{1}{2} {x^\omega_k}^T H^\omega_k x^\omega_k - (f^\omega_k)^T x^\omega_k\\
    \mathop{\textrm{s.t.}} \quad& A_{k} x^\omega_k + B_{k} z = 0 \;(y_k) \quad \forall k\in\mathcal{K},\label{eqn:lifted-constraint}
  \end{align}
\end{subequations}
where $z$ are the duplicates of coupling variables, $y$ are the dual variables for the coupling contraints, $A_k$ and $B_k$ are incidence matrices that extract coupling variables. ADMM iterations are performed by minimizing the augmented Lagrangian of \eqref{eqn:lifted} over $x$ and $z$ in turn and updating the dual $y$ based on the constarint violation of \eqref{eqn:lifted-constraint}. To make a fair comparison, both the subproblems of the overlapping scheme and the subproblems of ADMM are solved using {\tt IPOPT} \cite{wachter2006implementation} and modeled with {\tt JuMP} \cite{dunning2017jump}. Parallel synchronous schemes are implemented in {\tt Julia}.

The numerical results are shown in Fig. \ref{fig:benchmark}. As can be seen, the convergence rate (in terms of both the number of iteration and CPU time) of the overlapping scheme outperforms that of ADMM with the best penalty parameter. 
The superior convergence rate of the overlapping scheme is enabled by {\em exploiting the graph-structure} of the problem. In our view, ADMM exploits the structure in a more superficial manner (by partitioning the problem in the form of \eqref{eqn:lifted}).  In any case, the numerical results present strong evidence that the proposed overlapping scheme can outperform ADMM.

\begin{rmk}
  The claimed performance of the overlapping scheme only applies to graph-structured problems (Definition \ref{defn:bandwidth}(b)). Many problems in learning (e.g., linear regression problems of the form $\min_x \sum_i f_i(x)$ with no sparse coupling between $f_i(\cdot)$ and $x$) are not graph-structured. The graph structure of such problems has a small diameter, so there is no opportunity to exploit it efficiently. For such problems, our algorithm is not suitable. However, ADMM can still be effective in such problems.
\end{rmk}

\begin{figure}[!t]
  \centering
  \includegraphics[width=.35\textwidth]{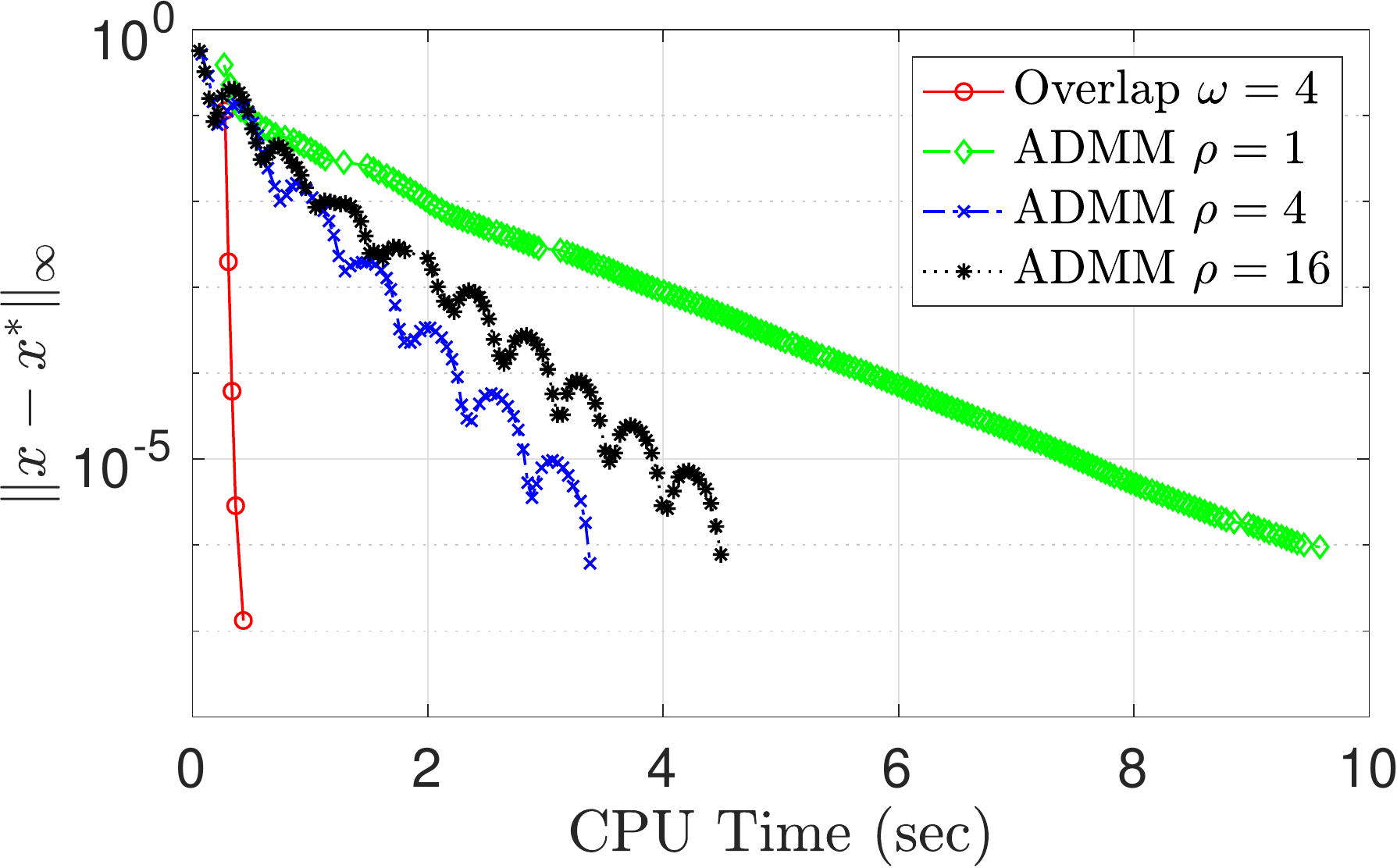}
  \caption{(Top) Error over CPU time ($c=0.1$, partition structure 1) for the proposed scheme with $\omega=4$ and ADMM with $\rho=1$, $4$, and $16$.}\label{fig:benchmark}
\end{figure}}

\section{Conclusions and Future Work}\label{sec:con}
We have presented a new decentralized computing paradigm that uses overlapping regions to promote convergence. This paradigm covers a spectrum of architectures that have fully centralized and decentralized architectures as extreme cases. We proved that the convergence rate of the proposed scheme improves exponentially with the size of the overlap and that convergence can be guaranteed under synchronous and asynchronous implementations. As part of future work, we are interested in exploring the convergence properties of the algorithm when applied to constrained optimization problems.

\section*{Acknowledgment}
We thank Guillaume Bal and Michel Schanen for comments and example cases. 
This material was based upon work
supported by the U.S. Department of Energy, Office of Science,
Office of Advanced Scientific Computing Research (ASCR), under
Contract DE-AC02-06CH11357 and partial NSF support
through award CNS-1545046 (MA) and EECS-1609183 (VZ).


\ifCLASSOPTIONcaptionsoff
\newpage
\fi



\appendix
\subsection{Proof of Lemma \ref{lem:bandwidth}}\label{prop:bandiwdth-pf}
We  observe that if $d_G(i,j)> \max\left(\mathcal{B}_G(A),\mathcal{B}_G(B)\right)$, then $(A)_{i,j}=(B)_{i,j}=0$, and thus $(A+B)_{i,j}=0$. This indicates that $\mathcal{B}_G(A+B)\leq \max\left(\mathcal{B}_G(A),\mathcal{B}_G(B)\right)$. Thus, \eqref{eqn:bandwidth-a} holds.

Suppose that $AB$ is diagonal. Then $\mathcal{B}_G(AB)=0$. Since $\mathcal{B}_G(A),\mathcal{B}_G(B)\geq 0$, \eqref{eqn:bandwidth-b} holds. Now we assume that $AB$ is not diagonal (i.e., $\mathcal{B}_G(AB)> 0$).  We know that there exists $(AB)_{i,j}$ such that $d_G(i,j)=\mathcal{B}_G(AB)>0$. Since $(AB)_{i,j} = \sum_{k=1}^n A_{ik} B_{kj}$, there exists $k$ such that $A_{ik}\neq 0$ and $B_{kj}\neq 0$. Furthermore, $d_G(i,k)\leq \mathcal{B}_G(A)$ and $d_G(k,j) \leq \mathcal{B}_G(B)$. Since $d_G(i,j)\leq d_G(i,k)+d_G(k,j)$, we have that $\mathcal{B}_G(AB) \leq \mathcal{B}_G(A) + \mathcal{B}_G(B)$. Thus, \eqref{eqn:bandwidth-b} holds.

\vspace{-0.1in}
{\subsection{Proof of Theorem \ref{thm:domdec}}\label{thm:domdec-pf}
We  observe that the following holds for $\lambda_{\text{mid}}:=(\lambda_{\min}+\lambda_{\max})/2$.
\begin{align}\label{eqn:domdec1}
  H^{-1} = \frac{1}{\lambda_{\text{mid}}}\left(I - (I - \frac{1}{\lambda_{\text{mid}}}H) \right)^{-1}
\end{align}
We also observe:
\begin{align}\label{eqn:sradi}
  \rho(1-\frac{1}{\lambda_{\text{mid}}}H) \leq \frac{\lambda_{\max}-\lambda_{\min}}{\lambda_{\max}+\lambda_{\min}}<1.
\end{align}
By applying \cite[Theorem 5.6.9 and Corollay 5.6.16]{horn1990matrix} to \eqref{eqn:domdec1}, we obtain:
\begin{align}\label{eqn:domdec2}
  H^{-1} = \frac{1}{\lambda_{\text{mid}}}\sum_{m=0}^\infty (I-\frac{1}{\lambda_{\text{mid}}}H)^m.
\end{align}
By Lemma \ref{lem:bandwidth}, $\mathcal{B}_G((I-\frac{1}{\lambda_{\text{mid}}} H)^m)\leq m\mathcal{B}_G(H)$. This indicates that $\left((I-\frac{1}{\lambda_{\text{mid}}}H)^m\right)_{i,j}=0$ if $d_G(i,j)>m\mathcal{B}_G(H)$. Thus, we have 
\begin{align}\label{eqn:domdec3}
  (H^{-1})_{i,j} &= \frac{1}{\lambda_{\text{mid}}}\sum_{m=\lceil\frac{ d_G(i,j)}{\mathcal{B}_G(H)}\rceil}^\infty \left((I-\frac{1}{\lambda_{\text{mid}}}H)^m\right)_{i,j};
\end{align}
and by applying the triangle inequality,
\begin{align}\label{eqn:domdec4}
  \left|(H^{-1})_{i,j}\right| &\leq \frac{1}{\lambda_{\text{mid}}}\sum_{m=\lceil\frac{ d_G(i,j)}{\mathcal{B}_G(H)}\rceil}^\infty \left|\left((I-\frac{1}{\lambda_{\text{mid}}}H)^m\right)_{i,j}\right|.
\end{align}
By Lemma \ref{lem:contraction} (see below) and \eqref{eqn:sradi}, we have:
\begin{subequations}\label{eqn:domdec5}
  \begin{align}
    \left|\left((I-\frac{1}{\lambda_{\text{mid}}}H)^m\right)_{i,j}\right| &\leq \rho\left(1-\frac{1}{\lambda_{\text{mid}}}H\right)^m\\
    &\leq \left(\frac{\lambda_{\max}-\lambda_{\min}}{\lambda_{\max}+\lambda_{\min}}\right)^m.
  \end{align}
\end{subequations}
From \eqref{eqn:domdec4}-\eqref{eqn:domdec5}, we obtain the following.
\begin{subequations}
  \begin{align}
    \left|\left(H^{-1}\right)_{i,j}\right| &\leq \frac{1}{\lambda_{\text{mid}}}\sum_{m=\lceil\frac{ d_G(i,j)}{\mathcal{B}_G(H)}\rceil}^\infty  \left(\frac{\lambda_{\max}-\lambda_{\min}}{\lambda_{\max}+\lambda_{\min}}\right)^m\\
    &\leq \frac{1}{\lambda_{\min}}\left(\frac{\lambda_{\max}-\lambda_{\min}}{\lambda_{\max}+\lambda_{\min}}\right)^{\frac{ d_G(i,j)}{\mathcal{B}_G(H)}}
  \end{align}
\end{subequations}

\begin{lem}\label{lem:contraction}
  Let $X\in\mathbb{R}^{n\times n}$ be symmetric. Then we have:
  \begin{align}
    \vert( X^k )_{ij} \vert\leq \rho(X)^k.
  \end{align}
\end{lem}
\begin{IEEEproof}
  There exist a diagonal matrix $\Lambda\in\mathbb{R}^{n\times n}$ and an orthogonal matrix $P\in\mathbb{R}^{n\times n}$ such that $X = P \Lambda P^T$ (see \cite[Theorem 4.1.5]{horn1990matrix}). Furthermore, $X^k = P^T \Lambda^k P$, and 
  \begin{subequations}
    \begin{align}
      |(X^k)_{i,j}| &= |\sum_{l=1}^n \lambda^k_l (P)_{i,l} (P)_{j,l}|\\
      &\leq \rho(X)^k\sum_{l=1}^n  |(P)_{i,l}|| (P)_{j,l}|
    \end{align}
  \end{subequations}
  We know that each row of $P$ is a unit vector on $\mathbb{R}^{n}$. Thus, $\sum_{l=1}^n  |(P)_{i,l}|| (P)_{j,l}|\leq 1$ since it takes the form of an inner product of unit vectors. Therefore, $|(X^k)_{i,j}|\leq\rho(X)^k$.
\end{IEEEproof}}
\vspace{-0.1in}
\subsection{Proof of Lemma \ref{lem:conv}}\label{lem:conv-pf}
It is well known (see \cite[Proposition 6.1, p144]{bertsekas1989parallel}) that $\rho(S^\omega)<1$ if and only if $I-S^\omega$ is nonsingular and the iteration of the form \eqref{eqn:jacobi} converges to its fixed point. {This follows from the fact that the power of the Jordan canonical form of $S^\omega$ converges to zero if and only if $\rho(S^\omega)<1$.} That is, the iteration converges to the solution of 
\begin{align}\label{eqn:jacobi-fixed}
  x = S^\omega x + U^\omega f.
\end{align}
By the assumption that $H$ is PD, \eqref{eqn:lin} has a unique solution. For both conditions \ref{lem:conv-a} and \ref{lem:conv-b} of Lemma \ref{lem:conv}, $I-S^\omega$ is nonsingular, and thus \eqref{eqn:jacobi-fixed} has a unique solution. Hence, we need only to show that the solution of \eqref{eqn:lin} is the solution of \eqref{eqn:jacobi-fixed}.

Let $x^*$ be the solution of \eqref{eqn:lin}. Then we have 
\begin{align}
  H^\omega_{k} \{x^*\}_{V^\omega_k} = -H^\omega_{-k} \{x^*\}_{V\setminus V^\omega_k} + f^\omega_{k}
\end{align}
for $k\in\mathcal{K}$. Since $H^\omega_{k}$ is nonsingular, we also have 
\begin{align}\label{eqn:conv1}
  \{x^*\}_{V_k} =\left\{(H^\omega_{k})^{-1}\left(- H^\omega_{-k} \{x^*\}_{V \setminus V^\omega_k} + f_k^\omega\right)\right\}_{i\in V_k}.
\end{align}
Equation \eqref{eqn:conv1} takes the same form with the fixed-point \eqref{eqn:jacobi-indwise}. Furthermore, one can show that \eqref{eqn:conv1} for $k\in\mathcal{K}$ is equivalent to $x^*=S^\omega x^* + U^\omega f$. Therefore, $x^*$ solves \eqref{eqn:jacobi-fixed}. This proves that the unique solutions of \eqref{eqn:lin} and \eqref{eqn:jacobi-fixed} are equal.

\vspace{-0.1in}
\subsection{Proof of Lemma \ref{lem:decen}}\label{prop:decen-pf}
From the block structure of $S^\omega$ we obtain:
\begin{align}
  \Vert S^\omega \Vert_\infty = \max_{k\in\mathcal{K}}\max_{i\in V_k} \sum_{j=1}^n |(S^\omega_k)_{i,j} |.
\end{align}
Furthermore, by the definition of $S^\omega_k$, we have 
\begin{align}
  (S_k^\omega)_{i,j}=
  \begin{cases}
    0 & \text{if}\quad j\in V^\omega_k \\
    \left((H^\omega_{k})^{-1}H^\omega_{-k}\right)_{i,j} & \text{otherwise}.
  \end{cases}
\end{align}
Thus,
\begin{align}\label{eqn:ineq1}
  \sum_{j=1}^n |(S_k^\omega)_{i,j}| \leq \sum_{j\in V\setminus V^\omega_k} \left|\left((H^\omega_{k})^{-1}H^\omega_{-k}\right)_{i,j}\right|.
\end{align}
If $V^\omega_k = V$, we have that $\sum_{j=1}^n |(S_k^\omega)_{i,j}| = 0$. Now suppose that $V^\omega_k\neq V$. Then we have the following.
\begin{subequations}
  \begin{align}
    \left|\left((H^\omega_{k})^{-1} H^\omega_{-k}\right)_{i,j}\right|
    &= \left|\sum_{l\in  V^\omega_k} \left((H^\omega_{k})^{-1}\right)_{i,l}\left(H^\omega_{-k}\right)_{l,j}\right|\\
    &\leq \sum_{l\in  V^\omega_k} \left|\left((H^\omega_{k})^{-1}\right)_{i,l}\right|\left|\left(H^\omega_{-k}\right)_{l,j}\right|\label{eqn:decen-sum}
  \end{align}
\end{subequations}
The eigenvalues of $H^\omega_{k}$ are on $[\lambda^k_{\min},\lambda^k_{\max}]$, and we have that $\lambda^k_{\min}>0$ since $H$ is PD. Thus, by applying Theorem \ref{thm:domdec}:
\begin{align}
  \left|\left((H^\omega_{k})^{-1}\right)_{i,l}\right|
  \leq \frac{1}{\lambda^k_{\min}} \left(\frac{\lambda^k_{\max}-\lambda^k_{\min}}{\lambda^k_{\max}+\lambda^k_{\min}}\right)^{d_G(i,l)/\mathcal{B}_G(H^\omega_{k})}.
\end{align}
If $(H^\omega_{-k})_{l,j}\neq 0$, then $d_G(l,j)\leq \mathcal{B}_G(H^\omega_{k})$. Thus, in the summation on \eqref{eqn:decen-sum}, we may  consider only  $l\in V^\omega_k$ such that $d_G(l,j)\leq \mathcal{B}_G(H^\omega_{k})$. This yields $d_G(i,j)\leq d_G(i,l)+\mathcal{B}_G(H^\omega_{k})$. By noting that $i\in V_k$ and $j\in V\setminus V^\omega_k$
and using the definition of $V^\omega_k$ in \eqref{eqn:ovlblk}, we obtain $d_G(i,j) \geq \omega+1$. Therefore,
\begin{align}
  \left|\left((H^\omega_{k})^{-1}\right)_{i,l}\right|
  \leq \frac{1}{\lambda^k_{\min}} \left(\frac{\lambda^k_{\max}-\lambda^k_{\min}}{\lambda^k_{\max}+\lambda^k_{\min}}\right)^{\frac{\omega+1}{\mathcal{B}_G(H^\omega_{k})}-1}.
\end{align}
Using the above observation, we can establish that
\begin{align}\label{eqn:ineq2}
  &\left|\left((H^\omega_{k})^{-1} H^\omega_{-k}\right)_{i,j}\right| \\
  &\leq \sum_{l\in V^\omega_k}\left|(H^\omega_{-k})_{l,j}\right|\frac{1}{\lambda^k_{\min}} \left(\frac{\lambda^k_{\max}-\lambda^k_{\min}}{\lambda^k_{\max}+\lambda^k_{\min}}\right)^{\frac{\omega+1}{\mathcal{B}_G(H^\omega_{k})}-1}.\nonumber
\end{align}
By applying \eqref{eqn:ineq2} to \eqref{eqn:ineq1}, we obtain:
\begin{subequations}\label{eqn:nonempty}
  \begin{align}
    \sum_{j=1}^n |(S_k^\omega)_{i,j}|
    &\leq \left(\sum_{j\in V\setminus V^\omega_k}
    \sum_{l\in V^\omega_k}\left|(H^\omega_{-k})_{l,j}\right|\right)\\
    &\qquad \times \frac{1}{\lambda^k_{\min}} \left(\frac{\lambda^k_{\max}-\lambda^k_{\min}}{\lambda^k_{\max}+\lambda^k_{\min}}\right)^{\frac{\omega+1}{\mathcal{B}_G(H^\omega_{k})}-1}\nonumber\\
    &= \frac{R_k}{\lambda^k_{\min}} \left(\frac{\lambda^k_{\max}-\lambda^k_{\min}}{\lambda^k_{\max}+\lambda^k_{\min}}\right)^{\frac{\omega+1}{\mathcal{B}_G(H^\omega_{k})}-1}
  \end{align}
\end{subequations}
for $i\in V_k$. Thus,
\begin{align}\label{eqn:Sbound}
  \max_{i\in V_k} \sum_{j=1}^n |(S_k^\omega)_{i,j}|\leq \frac{R_k}{\lambda^k_{\min}} \left(\frac{\lambda^k_{\max}-\lambda^k_{\min}}{\lambda^k_{\max}+\lambda^k_{\min}}\right)^{\frac{\omega+1}{\mathcal{B}_G(H^\omega_{k})}-1}.
\end{align}
Note that \eqref{eqn:Sbound} holds regardless if $V=V^\omega_k$ or not.
Finally,
\begin{subequations}
\begin{align}
  \Vert S^\omega\Vert_\infty
  & \leq \max_{k\in\mathcal{K}}\max_{i\in V_k}\sum_{j=1}^n |(S_k^\omega)_{i,j}|\\
  & \leq \max_{k\in\mathcal{K}} \frac{R_k}{\lambda^k_{\min}} \left(\frac{\lambda^k_{\max}-\lambda^k_{\min}}{\lambda^k_{\max}+\lambda^k_{\min}}\right)^{\frac{\omega+1}{\mathcal{B}_G(H^\omega_{k})}-1}.
\end{align}
\end{subequations}

\subsection{Proof of Theorem \ref{thm:async}}\label{prop:async-pf}
For a given $x^{(0)}\in\mathbb{R}^n$, we define the following:
\begin{align}
  X(\ell) := \{x\in\mathbb{R}^n\mid \Vert x-x^*\Vert_\infty \leq \left(\Vert S^\omega\Vert_\infty\right)^\ell\Vert x^{(0)}-x^*\Vert_\infty \}
\end{align}
for $\ell=0,1,\cdots$, where $S^\omega$ is defined in \eqref{eqn:SUk}-\eqref{eqn:SU}. Furthermore, we can consider an equivalent representation $X(\ell)= X_1(\ell)\times\cdots\times X_K(\ell)$, where
\begin{align}
  \begin{split}
    &X_k(\ell)= \left\{x_k\in\mathbb{R}^{|V_k|}:\right.\\
    &\left.\Vert x_k-\{x_i^*\}_{i\in V_k}\Vert_\infty \leq \left(\Vert S^\omega\Vert_\infty\right)^\ell\Vert x^{(0)}-x^*\Vert_\infty \right\}.
  \end{split}
\end{align}
Since $(I-S^\omega)x^*=U^\omega f$, for any $y\in X(\ell)$, we have that
\begin{subequations}\label{eqn:jacobi-y-bound}
  \begin{align}
    \Vert S^\omega y+U^\omega f-x^*\Vert_\infty
    &=\Vert S^\omega  \left(y-x^*\right)\Vert_\infty\\
    &\leq \left(\Vert S^\omega \Vert_\infty\right)^{\ell+1} \Vert x^{(0)}-x^*\Vert_\infty.
  \end{align}
\end{subequations}
Thus, we observe that for any $y\in X(\ell)$, we have that $ S^\omega y+U^\omega f\in X(\ell+1)$. In other words, for any $k\in\mathcal{K}$,
\begin{align}\label{eqn:progress}
  S^\omega_k y + U^\omega_k f \in X_k(\ell+1).
\end{align}

Based on this observation, we show the following. {\em Claim:} For each $\ell=0,1,\cdots$, there exists $t_\ell$ such that if $t\geq t_\ell$, then $x^{(t)},x^{k,(t)}\in X(\ell)$ for $k\in\mathcal{K}$, $t\in\mathcal{T}_k$, and $t\geq t_\ell$. We use mathematical induction. First, if $\ell=0$, we can choose $t_0=0$. By \eqref{eqn:progress}, we may consider \eqref{eqn:async} as a mapping from $X(0)$ to $X(0)$. Thus, we have that $x^{(t)}\in X(0)$ for any $t$. Since $\{x_{k'}^{(\tau_{k,k'}(t))}\mid t\in\mathcal{T}_{k}\}$ is a subset of $\{x_{k'}^{(t)}\mid t\in \mathcal{T}\}$ for any $k,k'\in\mathcal{K}$, $x^{k,(t)}$ is also in $X(0)$ for any $k\in\mathcal{K}$ and $t\in\mathcal{T}_k$. This proves that the claim holds for $\ell=0$.

Now we assume that the claim holds for $\ell$. By the induction hypothesis, $x^{k,(t)}\in X(\ell)$ for $t\geq t_\ell$ and $t\in\mathcal{T}_k$. If we choose $t'_{\ell+1}=\max_{k\in\mathcal{K}} \min \left\{t\in\mathcal{T}_k\mid t\geq t_\ell\right\}$, by \eqref{eqn:progress} we have $x^{(t)}\in X(\ell+1)$ for $t\geq t'_{\ell+1}$. Furthermore, since $\tau_{k,k'}(t)\rightarrow\infty$ as $t\rightarrow\infty$, there exists $t_{\ell+1} $ such that if $t\geq t_{\ell+1}$, $\tau_{k,k'}(t)\geq t'_{\ell+1}$ for any $k,k'\in\mathcal{K}$. Thus, we can establish that for $t\geq t_{\ell+1}$, $x^{k,(t)}\in X(\ell+1)$ for any $k\in\mathcal{K}$. The induction is complete, and we establish that there exists $\{t_\ell\}_{\ell=0,1,\cdots}$ such that $x^{(t)}\in X(\ell)$ if $t\geq t_\ell$. {We consider a mapping $t\mapsto\ell$ such that $\ell = \max \{\ell: t_\ell \leq t\}$ holds. With this mapping, we can write $x^{(t)}\in X(\ell)$, and this implies} that 
 $ \Vert x^{(t)}- x^*\Vert_\infty \leq \left(\Vert S^\omega \Vert_\infty \right)^\ell\Vert x^{(0)}-x^*\Vert$. Since $\ell\rightarrow\infty$ as $t\rightarrow\infty$, we obtain $x^{(t)}\rightarrow x^*$ as $t\rightarrow\infty$. 
\vspace{-0.1in}
{
\subsection{Generalization of Theorem \ref{thm:domdec}}           \label{apx:generalize}
  Here we present a theorem that generalizes the results in Theorem \ref{thm:domdec}, where we do not assume the PD of $H$. This result suggests that the convergence analysis of the algorithm can also be derived from Theorem \ref{thm:domdec-gen} without a PD assumption. Such examples arise in game theory and dynamical systems. 

\begin{thm}\label{thm:domdec-gen}
  Consider a matrix $H\in\mathbb{R}^{n\times n}$ and a graph $G(V,E)$ with $\{1,2,\cdots,n\}\subseteq V$. Suppose that the eigenvalues of $H$ are in $\{\lambda\in\mathbb{C} \mid \vert \lambda - z \vert \leq R\}$ for some $R\in\mathbb{R}_{>0}$ and $z\in\mathbb{C}\setminus\{0\}$ such that $R<|z|$. For any given $\epsilon\in(0,1-R/|z|)$, there exists $C=C(H,\epsilon)\in\mathbb{R}_{>0}$ such that
  \begin{align}
    \left|\left(\left(I-\frac{1}{z} H \right)^m\right)_{i,j}\right| \leq C \left(\frac{R}{|z|}+\epsilon\right)^m \label{eqn:domdec-gen-a}
  \end{align}
  holds for any $i,j\in\{1,2,\cdots,n\}$ and $m=1,2,\cdots$. Furthermore, for such $C$ and $\epsilon$, the following holds for any $i,j\in\{1,2,\cdots,n\}$:
  \begin{align}
    &\left|(H^{-1})_{i,j}\right| \leq \frac{C}{(1-\epsilon)|z|-R} \left(\frac{R}{|z|}+\epsilon\right)^{d_G(i,j)/\mathcal{B}_G(H)}.     \label{eqn:domdec-gen-b}
  \end{align}
\end{thm}
\begin{IEEEproof}
Any eigenvalue $\lambda$ of $H$ satisfies $|\lambda-z| \leq R$ and  is equivalent to
  $\left|\lambda/z- 1 \right|\leq R/|z|$.
Thus, the eigenvalues of $I-({1}/{z})H$ lie on $\{\lambda\in\mathbb{C}: |\lambda|\leq {R}/{|z|}\}$. This yields that $\rho(I-(1/z)H)<1$. By applying \cite[Corollary 5.6.13]{horn1990matrix}, we can show that there exists $C$ such that \eqref{eqn:domdec-gen-a} holds for given $\epsilon$.

We  observe that
\begin{align}\label{eqn:domdec-gen-1}
  H^{-1} = \frac{1}{z}\left(I - (I - \frac{1}{z}H) \right)^{-1},
\end{align}
and we know that $\rho((I - \frac{1}{z}H))<1$. By applying \cite[Theorem 5.6.9 and Corollay 5.6.16]{horn1990matrix}, we have that
\begin{align}\label{eqn:domdec-gen-2}
  H^{-1} = \frac{1}{z}\sum_{m=0}^\infty (I-\frac{1}{z}H)^m.
\end{align}
By Lemma \ref{lem:bandwidth}, $\mathcal{B}_G((I-\frac{1}{z} H)^m)\leq m\mathcal{B}_G(H)$. This indicates that $\left((I-\frac{1}{z}H)^m\right)_{i,j}=0$ if $d_G(i,j)>m\mathcal{B}_G(H)$. Thus, 
\begin{align}\label{eqn:domdec-gen-3}
  (H^{-1})_{i,j} &= \frac{1}{z}\sum_{m=\lceil\frac{ d_G(i,j)}{\mathcal{B}_G(H)}\rceil}^\infty \left((I-\frac{1}{z}H)^m\right)_{i,j};
\end{align}
and by applying the triangle inequality,
\begin{align}\label{eqn:domdec-gen-4}
  \left|(H^{-1})_{i,j}\right| &\leq \frac{1}{|z|}\sum_{m=\lceil\frac{ d_G(i,j)}{\mathcal{B}_G(H)}\rceil}^\infty \left|\left((I-\frac{1}{z}H)^m\right)_{i,j}\right|.
\end{align}
Using \eqref{eqn:domdec-gen-a}, we obtain
\begin{subequations}
  \begin{align}
    \left|\left(H^{-1}\right)_{i,j}\right| &\leq \frac{C}{|z|} \sum_{m=\lceil\frac{ d_G(i,j)}{\mathcal{B}_G(H)}\rceil}^\infty \left(\frac{R}{|z|}+\epsilon\right)^{m}\\
    &\leq \frac{C}{(1-\epsilon)|z|-R}\left(\frac{R}{|z|}+\epsilon\right)^{\frac{ d_G(i,j)}{\mathcal{B}_G(H)}}.
  \end{align}
\end{subequations}
\end{IEEEproof}}
\vspace{-0.1in}

\bibliographystyle{IEEEtran}
\bibliography{tcns}




%

\vspace{-0.3in}
\begin{IEEEbiographynophoto}{Sungho Shin}
  received his B.S. in chemical engineering and mathematics from Seoul National University, South Korea, in 2016. He is currently a Ph.D. candidate in the Department of Chemical and Biological Engineering at the University of Wisconsin-Madison. His research interests include control theory, optimization algorithms, and complex networks.
\end{IEEEbiographynophoto}
\vspace{-0.3in}
\begin{IEEEbiographynophoto}{Victor M. Zavala}
  is the Baldovin-DaPra Associate Professor in the Department of Chemical and Biological Engineering at the University of Wisconsin-Madison. He holds a B.Sc. degree from Universidad Iberoamericana and a Ph.D. degree from Carnegie Mellon University, both in chemical engineering. He is an associate editor for the Journal of Process Control and a technical editor of Mathematical Programming Computation. His research interests are in the areas of energy systems, high-performance computing, stochastic programming, and predictive control.
\end{IEEEbiographynophoto}
\vspace{-0.3in}
\begin{IEEEbiographynophoto}{Mihai Anitescu}
  is a senior computational mathematician in the Mathematics and Computer Science Division at Argonne National Laboratory and a professor in the Department of Statistics at the University of Chicago. He obtained his engineer diploma (electrical engineering) from the Polytechnic University of Bucharest in 1992 and his Ph.D. in applied mathematical and computational sciences from the University of Iowa in 1997. He specializes in the areas of numerical optimization, computational science, numerical analysis and uncertainty quantification. He is on the editorial board of Mathematical Programming A and B, SIAM Journal on Optimization, SIAM Journal on Scientiﬁc Computing, and SIAM/ASA Journal in Uncertainty Quantification; and he is a senior editor for Optimization Methods and Software.
\end{IEEEbiographynophoto}







\end{document}